\documentclass[10pt]{article}
\usepackage{amsfonts}
\usepackage{amscd}
\usepackage{amsmath}
\usepackage{latexsym, euscript}
\usepackage[latin1]{inputenc} 
\usepackage{setspace}

\def\d{\partial}

\numberwithin{equation}{section}

\newtheorem{thm}{Theorem}[section]
\newtheorem{cor}[thm]{Corollary}
\newtheorem{lem}[thm]{Lemma}
\newtheorem{prop}[thm]{Proposition}

\newtheorem{rem}[thm]{Remark}

\begin{document}

\newcommand{\CLIP}[2]{\left(#1,#2 \right)_{L^2}} 
\newcommand{\DLIP}[2]{\left(#1,#2 \right)_{L^2_h}} 
\newcommand{\DSIP}[2]{\left(#1,#2 \right)_{S_h}} 
\newcommand{\CLN}[1]{\left\| #1 \right\|_{L^2}}  
\newcommand{\CLNtx}[1]{\left\| #1 \right\|_{L^2_t L^2_x}}  
\newcommand{\CLNStx}[1]{\left\| #1 \right\|_{L^2_t L^2_x}^2}  
\newcommand{\CLNxt}[1]{\left\| #1 \right\|_{L^2_x L^2_t}}  
\newcommand{\DLN}[1]{\left\| #1 \right\|_{L^2_h}}  
\newcommand{\DLNt}[1]{\left\| #1 \right\|_{L^2_k}}  
\newcommand{\DLNSxh}[1]{\left\| #1 \right\|_{L^2_{x,h}}^2}  
\newcommand{\DLNStk}[1]{\left\| #1 \right\|_{L^2_{t,k}}^2}  
\newcommand{\DLNtx}[1]{\left\| #1 \right\|_{L^2_{t,k} L^2_{x,h}}}  
\newcommand{\DLNxt}[1]{\left\| #1 \right\|_{L^2_{x,h} L^2_{t,k}}}  
\newcommand{\DSN}[1]{\left\| #1 \right\|_{S_h}}  

\newcommand{\CLNS}[1]{\left\| #1 \right\|_{L^2}^2}  
\newcommand{\CLNSx}[1]{\left\| #1 \right\|_{L^2_x}^2}  
\newcommand{\CLNSt}[1]{\left\| #1 \right\|_{L^2_t}^2}  
\newcommand{\DLNS}[1]{\left\| #1 \right\|_{L^2_h}^2}  
\newcommand{\DSNS}[1]{\left\| #1 \right\|_{S_h}^2}  
\newcommand{\SupNorm}[1]{\left\| #1 \right\|_{\infty}}  
\newcommand{\SupNormtx}[1]{\left\| #1 \right\|_{L^{\infty}_t L^{\infty}_x}}  
\newcommand{\SupNormx}[1]{\left\| #1 \right\|_{ L^{\infty}_x}}  
\newcommand{\SupNormt}[1]{\left\| #1 \right\|_{L^{\infty}_t}}  
\newcommand{\SupNormxt}[1]{\left\| #1 \right\|_{L^{\infty}_x L^{\infty}_t}}  
\newcommand{\SupNormtLNx}[1]{\left\| #1 \right\|_{L^{\infty}_t L^2_x}}
\newcommand{\LNxSupNormt}[1]{\left\| #1 \right\|_{L^2_x L^{\infty}_t}}

\newcommand{\IHN}[1]{\left \langle #1 \right \rangle}  



\title{\bf Solutions of the Nonlinear Schr\"odinger Equation with Prescribed Asymptotics at Infinity} 
\author
{ 
\\ 
{\sc John Gonzalez}
\thanks{Supported in part by NSF grant DMS-0901443} 
\\
Department of Mathematics, Northeastern University\\
Boston, MA 02115 
}
\date{} 
\maketitle 
 
\begin{abstract}
\noindent We prove local existence and uniqueness of solutions for the one-dimensional nonlinear Schr\"odinger (NLS) equations $iu_t + u_{xx} \pm |u|^2 u = 0$ in classes of smooth functions that admit an asymptotic expansion at infinity in decreasing powers of $x$.  We show that an asymptotic solution differs from a genuine solution by a Schwartz class function which solves a generalized version of the NLS equation.  The latter equation is solved by discretization methods.  The proofs closely follow previous work done by the author and others on the Korteweg-De Vries (KdV) equation and the modified KdV equations.  
\end{abstract}

\section{Introduction}

In this article we consider the one-dimensional focusing and defocusing nonlinear Schr\"odinger equations 

\begin{eqnarray}
 \left\{ 
\begin{array}{rcl}
iw_t + w_{xx} + \mu |w|^2w &= &0 \label{NLS} \\
w_{|t=0}& =& w_0(x) 
\end{array} 
\right. 
\end{eqnarray}

\noindent where $\mu = \pm 1$ and the initial data $w_0$ has prescribed asymptotic expansions at $+\infty$ and/or $-\infty$.  Specifically we shall construct local (in time) solutions to (\ref{NLS}) that lie in the spaces $S^{\beta}\left(\mathbb{R}\times I \to \mathbb{C}\right)$ which are complex analogues to those spaces defined in \cite{Kap-Per-Shub-Top} by T. Kappeler, P. Perry, M. Shubin, and P. Topalov and used also by the author in \cite{Gonz}.  These spaces are defined as follows:  \\
\indent Let $I \subset \mathbb{R}$ be an interval and $\beta \in \mathbb{R}$ be given.  Denote by $S^{\beta}\left(\mathbb{R}\times I \to \mathbb{C}\right)$ (or simply $S^{\beta}\left(\mathbb{R}\times I \right)$) the linear space of $C^{\infty}\left(\mathbb{R}\times I \to \mathbb{C}\right)$ functions having asymptotic expansions at $\pm\infty$ given by $f(x,t) \sim \sum_{k=0} ^{\infty} \left( a_k ^{+}(t) + i b_k^{+}(t) \right) x^{\beta_k}$ as $x \to + \infty$ and $f(x,t) \sim \sum_{k=0} ^{\infty} \left( a_k ^{-}(t)+ i b_k^{-}(t) \right)x^{\beta_k}$ as $x \to -\infty$ where $a_k^{\pm}, b_k^{\pm}\in C^{\infty}(I \to \mathbb{R})$ and $\beta = \beta_0 > \beta_1 > \cdots$ with $\lim_{k \to \infty} \beta_k = -\infty$.  By definiton, the asymptotic relation $\sim$ means that for every compact interval $J\subset I$, and integers $N,i,j \geq 0$ there exists $C_{J,N,i,j} > 0$ such that for any $ \pm x \geq 1$ and $t\in J$ we have.  \\
$$\left| \partial^i_t \partial^j_x\left( f(x,t) - \sum_{k=0} ^N \left( a_k ^{\pm}(t)+i b_k^{\pm}(t) \right) (\pm x)^{\beta_k}\right) \right| \leq C_{J,N,i,j} \left| x \right| ^{\beta_{N+1} - j}$$
\noindent  We denote by $S^{-\infty}\left(\mathbb{R}\times I \to \mathbb{C} \right)$ (or simply $S^{-\infty}\left(\mathbb{R}\times I \right)$) the space of $C^{\infty}\left(\mathbb{R}\times I \to \mathbb{C}\right)$ functions having asymptotic expansions at $\pm \infty$ which are identically zero.  Analogously, we define the spaces $S^{\beta}\left(\mathbb{R} \to \mathbb{C} \right)$ and $S^{-\infty}\left(\mathbb{R}\to \mathbb{C}\right)$ as the space of functions $f(x) \in C^{\infty}(\mathbb{R} \to \mathbb{C})$ having such asymptotic expansions where the coefficients $a_k^{\pm}, b_k^{\pm}$ are constants independent of $t$.  We shall construct solutions $w(x,t) \in S^{\beta}(\mathbb{R}\times I)$ for (\ref{NLS}) with initial data $w_0\in S^{\beta}(\mathbb{R})$ when $\beta \leq 0$. \\
\indent  If $w(x,t) \in S^{\beta}(\mathbb{R}\times I)$ is a solution for (\ref{NLS}) then one expects its asymptotic expansions $\sum_{k=0} ^{\infty} \left( a_k ^{\pm}(t) +ib_k^{\pm}(t)\right)(\pm x)^{\beta_k}$, although not generally convergent, to give formal solutions (see lemma \ref{lem5.1}).  We define a pair of formal power series $\sum_{k=0} ^{\infty} \left( a_k ^{\pm}(t) +ib_k^{\pm}(t)\right) x^{\beta_k}$ to be a formal solution to (\ref{NLS}) if $\sum_{k=0} ^{\infty} \left( a_k ^{+}(t) +ib_k^{+}(t)\right) x^{\beta_k}$ and $\sum_{k=0} ^{\infty}\left( a_k ^{-}(t) +ib_k^{-}(t)\right) (-x)^{\beta_k}$ satisfy (\ref{NLS}) for all $t\in I$ when $x$ is taken as a formal variable and differentiation in $x$ is carried out in the ordinary way.  \\
\indent When $\beta > 0$ one can easily see that there are no such formal solutions to (\ref{NLS}) and hence no solutions in $S^{\beta}(\mathbb{R}\times I)$.  Indeed, if $\sum_{k=0} ^{\infty} \left( a_k ^{+}(t) +ib_k^{+}(t)\right) x^{\beta_k}$ satisfies (\ref{NLS}) formally where $\beta>0$ and $a_0^+ + i b_0^+ \ne 0$ then 
\begin{eqnarray*}
i \sum_{j=0}^{\infty}\left( \dot{a}_j^+(t) + i \dot{b}_j^+(t) \right)x^{\beta_j}+ \Big(\sum_{j=0}^{\infty} \left( a_j ^{+}(t) +ib_j^{+}(t)\right) \cdot \beta_j \cdot( \beta_j -1) x^{\beta_j-2}\Big) = \\
\mu \Big(\sum_{j=0}^{\infty} \left( a_j ^{+}(t) +ib_j^{+}(t)\right) x^{\beta_j}\Big)^2 \cdot \Big( \sum_{j=0}^{\infty} \left( a_j ^{+}(t) - ib_j^{+}(t)\right) x^{\beta_j} \Big)
\end{eqnarray*}
\noindent The largest exponent on the left side is $\beta_0$ and the largest exponent on the right side is $3\beta_0$ which is larger than $\beta_0$.  Therefore by equating the coefficients of $x^{3\beta_0 }$ one deduces that $0 = \mu (a_0^+ + i b_0^+)^2 (a_0^+ - i b_0^+)$ which implies that $a_0^+ = 0 = b_0^+$, a contradiction.  On the other hand when $\beta \leq 0$ there do exist such formal solutions defined for $t\in \mathbb{R}$ (see lemma \ref{lemA.2}) and therefore one can hope to find solutions in $S^{\beta}(\mathbb{R}\times I)$ for some interval $I$.  \\
\indent  For an arbitrarily chosen pair of such formal power series $\sum_{k=0} ^{\infty} \left( a_k ^{\pm}(t)+i b_k^{\pm}(t) \right) (\pm x)^{\beta_k}$ defined for $t\in I$ there exists a function $f(x,t)\in C^{\infty}(\mathbb{R}\times I \to \mathbb{C})$ asymptotic to the pair (see for example \cite{Shub} proposition 3.5).  The function $f$ is not unique but if $\sum_{k=0} ^{\infty} \left( a_k ^{\pm}(t)+i b_k^{\pm}(t) \right) (\pm x)^{\beta_k}$ is a formal solution then any such $f$ will be an asymptotic solution for (\ref{NLS}) (see lemma \ref{lemA.2}).  By definition an asymptotic solution is a function $f\in S^{\beta}(\mathbb{R}\times I \to \mathbb{C})$ such that \\
\[
\left\{ 
\begin{array}{rcl}
if_t + f_{xx} + \mu |f|^2f & \in & S^{-\infty} \left(\mathbb{R} \times I \to \mathbb{C}\right) \\
f_{|t=0} - w_0 & \in & S^{-\infty}\left(\mathbb{R} \to \mathbb{C}\right) 
\end{array} 
\right. 
\]  
\indent  Given an asymptotic solution $f(x,t) \in S^{\beta}(\mathbb{R}\times I)$ for (\ref{NLS}) one can attempt to construct a genuine solution $w(x,t) \in S^{\beta}(\mathbb{R}\times I)$ to (\ref{NLS})  by constructing $u(x,t) \in S^{-\infty} \left(\mathbb{R} \times I\right)$ such that $w := f+u$ is a genuine solution of (\ref{NLS}).  If $w$ satisfies (\ref{NLS}) then $u$ must satisfy\\

\begin{eqnarray}
 \left\{ 
\begin{array}{rcl}
iu_t + u_{xx} + \mu \left( u^2 \bar{u} + u^2 \bar{f} +f^2\bar{u} + 2 u f\bar{u} + 2 u f \bar{f}\right) + g &= &0 \label{gNLS}\\
u_{|t=0} &=& u_0(x) 
\end{array} 
\right. 
\end{eqnarray}

\noindent where $u_0 = w_0-f(x,0) \in S^{-\infty}\left(\mathbb{R} \to \mathbb{C} \right)$ and $g \in S^{-\infty}\left(\mathbb{R} \times I \to \mathbb{C}\right)$ is the result of plugging $f$ into (\ref{NLS}). \\
\indent We shall prove existence of finite time solutions $u(x,t) \in S^{-\infty}\left(\mathbb{R} \times \left[ 0,T \right] \right) $ to (\ref{gNLS}) by using a finite difference method as in \cite{Bond, Gonz, Menik}.  Moreover, uniqueness will also be proven so that we shall show the following theorem: \\
\begin{thm} \label{thm1.1}
Let $f\in C^{\infty}\left(\mathbb{R}\times \left[\left. 0,\infty \right. \right] \to \mathbb{C}\right)$ be a function satisfying the property that for every compact $J \subset \left[\left. 0,\infty \right. \right]$ we have $\left| f^{(n)}(x,t)\right| = O\left(\left| x \right|^{-n} \right)$ uniformly for $t\in J$ and let $g$ be any function lying in $S^{-\infty}\left(\mathbb{R}\times\left[\left. 0,\infty \right. \right] \to \mathbb{C} \right)$.  Suppose $u_0 \in S^{-\infty}\left(\mathbb{R}\to \mathbb{C}\right)$.  Then there exists $T>0$ such that  (\ref{gNLS}) has a solution $u(x,t) \in S^{-\infty}(\mathbb{R}\times \left[0,T\right]\to \mathbb{C})$.  Moreover, the solution $u$ is unique in $S^{-\infty}(\mathbb{R}\times \left[0,T\right]\to \mathbb{C})$.  
\end{thm}

\indent The finite-time existence and uniqueness theorem for (\ref{gNLS}) will enable us to prove finite-time existence and uniqueness for (\ref{NLS}) in the space $S^{\beta}(\mathbb{R}\times \left[0,T\right])$ for $\beta \leq 0$ which can be stated as the following main theorem:  \\

\begin{thm} \label{thm1.2}
For any $\beta \leq 0$ and for any initial condition $w_0 \in S^{\beta}(\mathbb{R} \to \mathbb{C})$ there exists a $T>0$ and a unique solution $w(x,t) \in S^{\beta}(\mathbb{R}\times \left[0,T\right] \to \mathbb{C})$ of the initial value problem (\ref{NLS}).  Moreover, if $w_0 \sim \sum_{k=0}^{\infty} \left( a_k^{\pm} +i b_k^{\pm}\right)x^{\beta_k}$ and $j$ is the smallest index such that $a_j^{+} + i b_j^{+} \ne 0$ (resp. $a_j^{-}+i b_j^{-} \ne 0$) then the coefficient $a_j^{+}(t) + i b_j^{+}(t)$ (resp. $a_j^{-}(t)+i b_j^{-}(t)$) in the asymptotic expansion of the solution is a nonvanishing continuous function of $t$ and all preceeding coefficients are identically zero.
\end{thm}

\indent The second statement in theorem \ref{thm1.2} indicates that the asymptotic decay rate of the solution is determined throughout its time of existence by the leading exponents in the asymptotic expansion of its initial data.  In particular if $\beta_0 =0$ and $a_0^{\pm}\ne 0$ then the solution $w(x,t)$ for (\ref{NLS}) is asymptotically constant in $x$.  \\

\noindent  \textit{Related Work} \hspace{1mm} The results and methods presented here are most closely related to those presented by the author's article on the modified Korteweg-De Vries (mKdV) equation \cite{Gonz} and the precursor articles on the KdV equation by Bondareva, Shubin, and Menikoff \cite{Bond, Bond-Shub, Menik}.  In \cite{Gonz} local existence and uniqueness is proved for the mKdV equation in spaces $S^{\beta}\left(\mathbb{R} \times I \to \mathbb{R} \right)$ for $\beta \leq \frac{1}{2}$ and in \cite{Bond, Bond-Shub} global existence and uniqueness is proved for the KdV equation in $S^{\beta}(\mathbb{R}\times \mathbb{R} \to \mathbb{R})$ for $\beta <1$.  In \cite{Menik} global existence for the KdV equation is proved in slightly different spaces of functions whose asymptotic growth is of order $|x|^{\beta}$ for $\beta \leq 1$.  For the defocusing modified KdV equation T. Kappeler, P. Perry, M. Shubin, and P. Topalov in \cite{Kap-Per-Shub-Top} proved global existence and uniqueness in $S^{\beta}(\mathbb{R}\times \mathbb{R} \to \mathbb{R})$ (as well as other spaces) for $\beta < \frac{1}{2}$.  For the defocusing NLS equation (i.e. $\mu=-1$ ) and the closely related Gross-Pitaevskii equation there are several articles where solutions which are asymptotically constant in $x$ are considered.  We mention here articles \cite{Ala-Car, Beth-Saut, Gal04, Gal08, Gal, Ger, Goub, Gus-Nak-Tsa06, Gus-Nak-Tsa07, Mar, Zhid} and some of the references therein.  \\
\indent In section two we introduce a discretization of (\ref{gNLS}), give some general lemmas, and prove finite time existence for the discrete equation.  Section three contains various estimates which are necessary in order to pass from discrete solutions to smooth solutions of (\ref{gNLS}).  In section four we show how to pass from discrete solutions to smooth solutions by using a smoothing operator $I_h$, introduced by Stummel in \cite{Stum}.  The existence statements of theorems \ref{thm1.1} and \ref{thm1.2} are proved in section four and the uniqueness results are proved in section five. 
\section{Discretization of the Generalized NLS Equation}

\subsection{Definitions and General Setup}
We shall first split equation (\ref{gNLS}) into its real and imaginary components and discretize those two equations.  We will also define some discrete functional spaces and give some simple properties of the definitions. \\
\indent  Suppose that $\beta \leq 0$, $u(x,t) = u_1(x,t)+i u_2(x,t)$, $f(x,t) = f_1(x,t) + i f_2(x,t)$, and $g(x,t) = g_1(x,t) + i g_2(x,t)$ where for $l=1,2$ we have $u_l, g_l \in S^{-\infty}\left(\mathbb{R}\times I \to \mathbb{R} \right)$  and $f_l \in S^{\beta}\left( \mathbb{R} \times I \to \mathbb{R}\right)$.  Then we may write the real and imaginary parts of equation (\ref{gNLS}) as

\begin{equation}
\d_t u_2 - \d_x^2 u_1 - \mu \left( u_1^3 +u_1 u_2^2 + 3u_1^2 f_1 + u_2^2 f_1 + 2 u_1 u_2 f_2 + 3 u_1 f_1^2 + u_1 f_2^2 + 2 f_1 f_2 u_2 \right) -g_1 =0 \label{rgNLS}
\end{equation} 
 and
\begin{equation}
\d_t u_1 + \d_x^2 u_2 + \mu \left(u_1^2u_2 + u_2^3+ 2u_1 u_2 f_1+u_1^2 f_2 + 3 u_2^2 f_2+2f_1 f_2 u_1+u_2 f_1^2 + 3 u_2 f_2^2\right) + g_2 = 0 \label{igNLS}
\end{equation}
Obviously solving the system (\ref{rgNLS}), (\ref{igNLS}) is equivalent to solving the equation (\ref{gNLS}) so we shall focus our efforts on the system.  \\
\indent For now let us fix two mesh size numbers $0 < h,k <1 $ and let us denote $x_n := nh$ and $t_j := jk$ for each $n,j \in \mathbb{Z}$.   We shall let $\mathbb{R}_h$ and $\mathbb{R}_k$ denote the (discrete) collection of real numbers of the form $x_n$ and $t_j$ respectively and we shall refer to those sets and the cartesian product $\mathbb{R}_h \times \mathbb{R}_k$ as meshes.  If $\rho$ is any function defined on a mesh and taking values in either $\mathbb{R}$ or $\mathbb{R}^2$ then we will refer to $\rho$ as a mesh function.  Obviously any function defined on a continuum $\mathbb{R}$ or $\mathbb{R} \times \mathbb{R}$ (which we may call continuum functions) can also be considered as a mesh function by restricting its domain to the mesh.  If $\rho$ is a mesh function on $\mathbb{R}_h \times \mathbb{R}_k$ then we will ease some notation by writing $\rho_{n,j} :=  \rho(x_n,t_j)$ and $\rho_j := \rho(\cdot, t_j)$.    \\
\indent We shall make use of two discrete derivative operators $D_+$, $D_-$ that "differentiate" mesh functions $\rho$ defined on $\mathbb{R}_h$ (and hence they can also differentiate continuum functions $\rho$ defined on $\mathbb{R}$).  The operators are given by  
\begin{eqnarray*}
D_+ \rho(x) = \frac{\rho(x+h) - \rho(x)}{h} \hspace{1cm}
D_- \rho(x) = \frac{\rho(x) - \rho(x-h)}{h} \hspace{1cm}
\end{eqnarray*} 
\noindent We will also sometimes use shifting operators $E$ and $E^{-1}$ given by 
\begin{eqnarray*}
(E\rho)(x)  =  \rho(x+h) \hspace{3cm} (E^{-1}\rho)(x) = \rho(x-h)
\end{eqnarray*}
\noindent These two operators will only act on the $x$ variable of our functions $\rho(x,t)$.  We will also use the operator $D_{t,+}\eta(t) = \frac{\eta(y+k) - \eta(t)}{k}$.  The following properties of $D_+$, $D_-$, $D_{t,+}$, $E$, and $E^{-1}$ are immediate consequences of their definition:

\begin{enumerate}
\item If $\rho = \rho(x,t)$ then the operators $D_+$, $D_-$, $D_{t,+}$ and $E$ all commute when acting on $\rho(x,t)$.
\item $D_+(\nu \cdot \rho)(x_n) = \nu(x_n) D_+ \rho(x_n) + (E \rho)(x_n)D_+\nu(x_n)$
\item $D_-(\nu \cdot \rho)(x_n) = \nu(x_n) D_- \rho(x_n) + (E^{-1}\rho)(x_n)D_-\nu(x_n) $
\item For any $n\in\mathbb{N}$ we have $$D_+^n (\rho \cdot \nu \cdot \xi) = \sum_{i_1+i_2 + i_3 = n}c_{i_1,i_2,i_3} (E^{i_2 + i_3}D_+^{i_1}\rho)\cdot(E^{i_3}D_+^{i_2}\nu)\cdot(D_+^{i_3}\xi) $$ for some constants $c_{i_1,i_2,i_3}\in\mathbb{N}$.
\item If a continuum function $\rho : \mathbb{R}\to\mathbb{R}$ is differentiable on $\mathbb{R}$ then for each $x_n \in \mathbb{R}_h$ there exists $x\in\mathbb{R}$ where $x_n \leq x \leq x_{n+1}$ such that we have $D_+ \rho(x_n) = \frac{d}{dx}\rho(x)$.
\end{enumerate}

\indent In order to solve the system (\ref{rgNLS}), (\ref{igNLS}) we shall consider the following system of difference schemes, which are discrete versions of (\ref{rgNLS}) and (\ref{igNLS})

\begin{align}
& D_{t,+}(u_2)_j - D_+D_-(u_1)_{j+1}- \mu \left[ (u_1)_j^2 (u_1)_{j+1} + (u_1)_j (u_2)_j (u_2)_{j+1} + (f_2)_j^2 (u_1)_j + 3(f_1)_j (u_1)_j (u_1)_{j+1} \right. \nonumber \\
 & \left.  + (f_1)_j (u_2)_j(u_2)_{j+1} + 2 (f_2)_j (u_1)_j(u_2)_{j+1}+ 3 (f_1)_j^2 (u_1)_j + 2(f_1)_j (f_2)_j(u_2)_j\right] - (g_1)_j = 0 \label{drgNLS}
\end{align}

\begin{align}
& D_{t,+}(u_1)_j +  D_+D_-(u_2)_{j+1} + \mu \left[ (u_1)_j(u_2)_j(u_1)_{j+1}+ (u_2)_j^2(u_2)_{j+1} + 2(u_2)_j(f_1)_j(u_1)_{j+1} \right. \nonumber \\
& \left. + (u_1)_j(f_2)_j(u_1)_{j+1}+3(u_2)_j(f_2)_j(u_2)_{j+1}+2(f_1)_j(f_2)_j(u_1)_j+(f_1)_j^2(u_2)_j+3(f_2)_j^2(u_2)_j \right] + (g_2)_j = 0 \label{digNLS}
\end{align}

\noindent In order to rewrite this system of difference schemes in a more convenient and concise form we introduce a linear operator $Q_j$ on pairs of mesh functions $\left(\rho_1(x_n),\rho_2(x_n)\right)$ given by 

\begin{align}
 Q_j (\rho_1,\rho_2) & :=\left(D_+D_- \rho_2 +\mu (u_1)_j(u_2)_j \rho_1 + \mu (u_2)_j^2\rho_2 + 2 \mu(u_2)_j(f_1)_j \rho_1 \right. \nonumber \\ & + \mu(u_1)_j(f_2)_j \rho_1 + 3 \mu (u_2)_j(f_2)_j\rho_2 , 
  - D_+D_-\rho_1-\mu (u_1)_j^2\rho_1 -\mu (u_1)_j(u_2)_j \rho_2 \nonumber \\ 
 & \left.-3\mu (f_1)_j (u_1)_j \rho_1 - \mu(f_1)_j (u_2)_j \rho_2 -2\mu(f_2)_j(u_1)_j\rho_2  \right) \hspace{2mm}\label{Q_j}
\end{align}

\noindent Then the system (\ref{drgNLS}), (\ref{digNLS}) can be written in a shorter form as 

\begin{eqnarray}
 \left(I+kQ_j\right) \left( (u_1)_{j+1}, (u_2)_{j+1} \right) = \left( (u_1)_{j}, (u_2)_{j} \right) - k\mu \left( (f_1)_j(f_2)_j(u_1)_j + (f_1)_j^2(u_2)_j + 3(f_2)_j^2(u_2)_j \hspace{2mm}, \right. \nonumber \\
  \left. -(f_2)_j^2(u_1)_j -3(f_1)_j^2(u_1)_j-2(f_1)_j(f_2)_j(u_2)_j\right) -k \left( (g_2)_j, -(g_1)_j \right) \hspace{3mm}\label{dgNLS}
\end{eqnarray}

\noindent where $g_j \in S^{-\infty}\left(\mathbb{R}\right)$ is considered as a mesh function.  The task behind solving (\ref{dgNLS}) then is to show that one can invert the operator $I+k Q_j$, at least for some finite amount of time.  The invertibility will be possible only in certain function spaces, therefore we will now introduce an appropriate space.  \\
\indent Suppose $u(x_n) = \left(u_1(x_n),u_2(x_n)\right)$ and $v(x_n)=\left(v_1(x_n),v_2(x_n) \right)$ are mesh functions where $u_l,v_l : \mathbb{R}_h \to \mathbb{R}$ for $l=1,2$.  We consider the discrete inner products given by
\begin{eqnarray*}
\DLIP{u_1}{v_1} &=& \sum_{-\infty} ^{\infty}u_1(x_n)v_1(x_n) h   \\
\DSIP{u_1}{v_1} &=& \DLIP{ \IHN{x}  u_1}{ \IHN{x} v_1} + \DLIP{ \IHN{x} D_+ u_1}{ \IHN{x} D_+ v_1} + \DLIP{D_+^2 u_1}{D_+^2 v_1}  \\
\DLIP{u}{v} &=& \DLIP{u_1}{v_1}+\DLIP{u_2}{v_2} \\
\DSIP{u}{v} &=& \DSIP{u_1}{v_1}+\DSIP{u_2}{v_2}
\end{eqnarray*} 

\noindent where $\IHN{x} = \sqrt{x^2+1}$, and we define the corresponding norms and Hilbert spaces,
$$\DLNS{u} = \DLIP{u}{u} \hspace{1cm} L^2_h = \left\{ u(x_n) \hspace{1mm} mesh \hspace{1mm} functions \hspace{1mm}on \hspace{1mm}\mathbb{R}_h: \DLN{u} < \infty \right\}  $$ 
$$\DSNS{u} = \DSIP{u}{u} \hspace{1cm} S_h = \left\{ u(x_n)  \hspace{1mm} mesh \hspace{1mm} functions \hspace{1mm} on \hspace{1mm} \mathbb{R}_h : \DSN{u} < \infty \right\} $$

\noindent  where we shall understand from the context whether $L^2_h$ and $S_h$ refers to single-valued or double-valued mesh functions.  From the definitions of the $L^2_h$ inner product and its norm we have the properties: 
\begin{enumerate}
\item $\DLN{E\rho}=\DLN{\rho}$
\item If $\rho, \nu, D_- \rho, D_+ \nu \in L^2_h$ then $\DLIP{D_+ \nu}{\rho} = - \DLIP{\nu}{D_-\rho} $
\item For any $j,k \in \mathbb{N}$ we have $\DLN{D_+^j D_-^k \rho} \leq C\DLN{\rho}$ where $C$ is a constant depending on $h$.
\item For any $l,N \in\mathbb{N}$ and $0 \leq j \leq N$ we have $\DLN{\left( D_+^j\IHN{x}^N \right) E^l \rho} \leq C \DLN{\IHN{x}^{N-j} \rho}$ where $C$ is independent of $h$.
\end{enumerate}

\indent The first and second properties follow from simply reindexing and/or rearranging terms in the summation. The third property just requires use of the triangle inequality on each summand of $D_+^j D_-^k u$. For the fourth property we write the definition of $D_+$ and $\IHN{x}$ and use the fact that $h\in \left[0,1\right]$.
\subsection{Preliminary Lemmas}
\noindent The following proposition will be used frequently and often without reference.  For the proof we refer the reader to \cite{Gonz} appendix B. 

\begin{prop} \label{prop2.1}
Let $N,n \in \mathbb{N}$, $T>0$, $0<h_1,k_1<1$ and let $g\in S^{-\infty}(\mathbb{R} \times \left[ \left.-c,\infty \right)\right.)$ for some $c>0$.  There exists $C_{N,n} > 0$ such that $\DLN{\IHN{x}^N D_+^n g_j} < C_{N,n}$ for each $0<h\leq h_1$, $0<k\leq k_1$, and $0\leq t_j \leq T$. 
\end{prop}

\indent Another simple but important fact that we will frequently use is the following:\\
\noindent If $f\in C^{\infty}\left(\mathbb{R}\times \left(-\infty,\infty \right) \to \mathbb{R} \right)$ satisfies the property that for every $n\in \mathbb{N}$ and for every compact interval $J \subset \left(-\infty,\infty \right)$ we have $\frac{d^n}{dx^n}f(x,t) = O\left(\left| x \right|^{-n} \right)$ uniformly for $t \in J$ then for each $t_j\in J$ and for all $x\in \mathbb{R}$ we have 
$$\left| \frac{d^n}{dx^n}f(x, t_j ) \right| \leq C \left|\IHN{x}^{-n}\right|$$ where $C>0$ is independent of $k$ and $j$ (but C might depend on $J$).  This statement follows directly from the definitions of $O$ and $\IHN{x}$. \\
\indent The discrete Sobolev inequalities stated below will allow us to prove that the operators $I+kQ_j$ for $j \in \mathbb{N}$ are bounded below and are thus invertible.  These inequalities are stated and proven in \cite{Gonz} appendix B.  \\
\indent  As a notational remark, from now on we will let $C$ denote a constant whose value might change between consecutive inequalities but the variables that it depends on will often be noted by its indices for example as $C_{n,j,h}$ means some constant depending on $n,j,$ and $h$.

\begin{lem} \label{lem2.2}
For every $n\in \mathbb{N} $ there exists $C_n>0$ such that for every $h>0$ and for every mesh function $u:\mathbb{R}_h\to \mathbb{R}$ we have
\begin{enumerate}
\item   $\DLN{ D^k_+ u } \leq C_n\left( \DLN{ u} + \DLN{ D^n_+u} \right) $  \hspace{1cm} for $0 \leq k \leq n$
\item   $\SupNorm{ D^k_+ u } \leq C_n \left( \DLN{ u} + \DLN{ D^n_+u} \right)$  \hspace{1cm} for $0 \leq k < n$ 
\end{enumerate}
\end{lem}

\begin{cor} \label{cor2.3}
For all $N,l,j \in \mathbb{N}$, there exists $C_{N,j,l} > 0$ such that for any $h \in (0,1)$ and for all mesh functions $u : \mathbb{R}_h \to \mathbb{R}$ we have
\begin{enumerate}
\item   $\DLN{ \IHN{x}^N D^j_+ u } \leq C_{N,j,l} \left( \DLN{ \IHN{x}^N u} + \DLN{\IHN{x}^N D^{j+l}_+u} \right) $ \hspace{1cm} for $l\geq 0$
\item  $\SupNorm{ \IHN{x}^N D^j_+ u } \leq C_{N,j,l} \left( \DLN{ \IHN{x}^N u} + \DLN{\IHN{x}^N D^{j+l}_+u} \right) $ \hspace{1cm} for $l\geq 1$
\end{enumerate}
\end{cor}

\indent The next lemma will allow us to prove that the solutions stay bounded for finite time with respect to the Schwartz semi-norms.  The proof can be found in \cite{Menik}.  

\begin{lem} \label{lem2.4}
Suppose $P,Q$ are $C^1\left(\mathbb{R}\right)$, nondecreasing, positive functions, $\Delta t >0$, and for each $j \in \mathbb{N}$ we have $t_j := j \Delta t$.  Let $\eta:\left[ 0, T_0 \right] \rightarrow \mathbb{R}$ be an arbitrary function satisfying \\
\begin{equation}
\frac{\eta_{j+1}-\eta_j}{\Delta t} \leq P(\eta_j) \eta_{j+1} + Q(\eta_j) \label{dgron}
\end{equation}
for each $t_j,t_{j+1} \in \left[0, T_0\right] $ where $\eta_j := \eta(t_j)$, and suppose that $\eta_0 \leq K$ for some $K >0$.  Then there exists $0<T \leq T_0$ and $L,\epsilon >0$ all three depending on $K, P, $ and $Q$ such that if $\Delta t < \epsilon$ then $\eta_j \leq L$ for each $j$ where $t_j \leq T$.  Moreover, if $P$ and $Q$ are constants then we may take $T=T_0$. 
\end{lem}

\subsection{Finite Time Existence for Discrete Generalized NLS Equation in $S_h$}
We will now prove finite time existence for (\ref{dgNLS}).  The following lemma is the key estimate for establishing invertibility of the operator $I+kQ_j$ in the space $S_h$.  \\

\begin{lem} \label{lem2.5}
Suppose $T>0$, $h,k \in \left(0,1 \right)$ and that for each $j$ where $t_j\in \left[0,T\right]$ we have a given mesh function $u_j : \mathbb{R}_h \to \mathbb{R}^2$.  Define the operators $Q_j$ as in (\ref{Q_j}). Then there exists $C>0$ depending only on $f$ and $T$ but not on $h$, $k$, $j$, or the mesh functions $u_j$ such that for any mesh function $u : \mathbb{R}_h \to \mathbb{R}^2$ the inequality \\
\begin{equation}
\DSIP{Q_j u}{u} \geq -C \DSNS{u } \left(1+ \DSNS{ u_j }\right) \label{qbound}
\end{equation}
holds for each $j$ where $0 \leq t_j \leq T$.\\
\end{lem}

\noindent \textbf{Proof of Lemma 2.5} \hspace{1mm} We expand the left side of (\ref{qbound}) by using the definition of $Q_j$ and $\DSIP{\cdot}{\cdot}$ and obtain a sum of terms which can be bounded by the right side of (\ref{qbound}).  We will now show how to bound those terms by the right side of (\ref{qbound}) for some appropriate constant $C$.  Upon adding all the inequalities we will obtain inequality (\ref{qbound}).  Throughout we will denote $u(x_n)=\left(u_1(x_n),u_2(x_n) \right)$ and $u_j =\left( (u_1)_j, (u_2)_j \right)$.

\noindent \underline{\textit{Estimate for Term}} $\DSIP{D_+D_-u_2}{u_1}-\DSIP{D_+D_-u_1}{u_2}$: 
\begin{align*}
\DSIP{D_+D_-u_2}{u_1}-\DSIP{D_+D_-u_1}{u_2} & =  \DLIP{D_+D_-u_2}{\IHN{x}^2 u_1} - \DLIP{D_+D_-u_1}{\IHN{x}^2 u_2}  \\
                                            &    + \DLIP{D_+^2D_- u_2}{\IHN{x}^2 D_+ u_1} -\DLIP{D_+^2D_-u_1}{\IHN{x}^2D_+u_2} \\
                                            &    +\DLIP{D_+^3D_-u_2}{D_+^2u_1}-\DLIP{D_+^3D_-u_1}{D_+^2u_2} \\
                                            & =  - \DLIP{D_+u_2}{D_+\IHN{x}^2 Eu_1} - \DLIP{D_+u_2}{\IHN{x}^2D_+u_1} \\
                                            &    + \DLIP{D_+u_1}{D_+\IHN{x}^2Eu_2} + \DLIP{D_+u_1}{\IHN{x}^2D_+u_2} \\
                                            &    - \DLIP{D_+^2u_2}{D_+\IHN{x}^2D_+u_1} - \DLIP{D_+^2u_2}{E\IHN{x}^2D_+^2u_1}\\
                                            &    + \DLIP{D_+^2u_1}{D_+\IHN{x}^2D_+u_2} + \DLIP{D_+^2u_1}{E\IHN{x}^2D_+^2u_2} \\
                                            &    - \DLIP{D_+^3u_2}{D_+^3u_1} + \DLIP{D_+^3u_1}{D_+^3u_2} \\
                                            & =  - \DLIP{D_+u_2}{D_+\IHN{x}^2 Eu_1} + \DLIP{D_+u_1}{D_+\IHN{x}^2Eu_2} \\
                                            &    - \DLIP{D_+^2u_2}{D_+\IHN{x}^2D_+u_1} + \DLIP{D_+^2u_1}{D_+\IHN{x}^2D_+u_2} \\
                                            & \geq  -C \left[ \DLN{D_+u_2}\DLN{\IHN{x}u_1} + \DLN{D_+u_1}\DLN{\IHN{x}u_2} \right. \\
                                            &       + \left. \DLN{D_+^2 u_2}\DLN{\IHN{x}D_+u_1} + \DLN{D_+^2u_1}\DLN{\IHN{x}D_+u_2} \right] \\
                                            & \geq  -C \DSN{u_1}\DSN{u_2} \geq -C \DSNS{u}
\end{align*}

\noindent \underline{\textit{Estimate for Term}} $\DSIP{3\mu(u_2)_j (f_2)_j u_2}{u_1}$: \\ 
By definition of the norm we have 
\begin{eqnarray}
\DSIP{3\mu(u_2)_j (f_2)_j u_2}{u_1} & = & \DLIP{3\mu\IHN{x}(u_2)_j (f_2)_j u_2}{\IHN{x}u_1} + \DLIP{3\mu\IHN{x}D_+\left( (u_2)_j (f_2)_j u_2\right)}{\IHN{x}D_+u_1} \nonumber \\
                                    &   & + \DLIP{3\mu D_+^2\left( (u_2)_j (f_2)_j u_2\right)}{D_+^2 u_1} \label{expan}
\end{eqnarray}
For the first term of (\ref{expan})we have the bound
\begin{eqnarray*}
\DLIP{3\mu\IHN{x}(u_2)_j (f_2)_j u_2}{\IHN{x}u_1} & \geq & -C\SupNorm{(f_2)_j}\SupNorm{(u_2)_j}\DLN{\IHN{x}u_2}\DLN{\IHN{x}u_1} \\
                                                  & \geq  &-C \DSN{(u_2)_j}\DSN{u_2}\DSN{u_1} \\
                                                  & \geq & -C \DSN{u_j}\DSNS{u} \geq -C \DSNS{u} \left( 1+\DSNS{u_j} \right)
\end{eqnarray*}
For the second term of (\ref{expan}) we use the product rule for $D_+$ to obtain
\begin{align*}
 & \DLIP{3\mu\IHN{x}D_+\left( (u_2)_j (f_2)_j u_2\right)}{\IHN{x}D_+u_1}  \hspace{1cm}\\
& =  \sum_{i_1+i_2+i_3=1} c_{i_1,i_2,i_3} \DLIP{E^{i_2+i_3} D_+^{i_1}(u_2)_j E^{i_3}D_+^{i_2}(f_2)_j D_+^{i_3}u_2}{\IHN{x}D_+u_1}
\end{align*}
For $i_1=1$ and other indices zero we have
\begin{eqnarray*}
3\mu c_1 \DLIP{\IHN{x}(f_2)_j u_2 D_+(u_2)_j}{\IHN{x}D_+u_1} & \geq & -C \SupNorm{(f_2)_j} \SupNorm{D_+(u_2)_j} \DLN{\IHN{x}u_2} \DLN{\IHN{x}D_+u_1} \\
                                                             & \geq & -C \DSN{(u_2)_j} \DSN{u_2}\DSN{u_1} \geq -C \DSN{u_j}\DSNS{u} \\
                                                             & \geq & -C \DSNS{u} \left(1+ \DSNS{u_j} \right)
\end{eqnarray*}
For $i_2=1$ and other indices zero we have
\begin{eqnarray*}
3\mu c_2 \DLIP{\IHN{x}u_2E(u_2)_jD_+(f_2)_j}{\IHN{x}D_+u_1} & \geq & -C \SupNorm{D_+(f_2)_j} \SupNorm{E(u_2)_j} \DLN{\IHN{x} u_2} \DLN{\IHN{x}D_+u} \\
                                                            & \geq & -C \SupNorm{\frac{d}{dx}(f_2)_j} \SupNorm{(u_2)_j} \DSN{u_2} \DSN{u_1} \\
                                                            & \geq & -C \DSN{(u_2)_j} \DSNS{u} \geq -C\DSNS{u} \left(1+\DSNS{u_j} \right)
\end{eqnarray*}
For $i_3=1$ and other indices zero we have
\begin{align*}
& 3\mu c_3 \DLIP{\IHN{x}E(u_2)_j E(f_2)_j D_+u_2}{\IHN{x}D_+u_1} \\
& \geq  -C\SupNorm{E(u_2)_j} \SupNorm{E(f_2)_j} \DLN{\IHN{x}D_+u_2} \DLN{\IHN{x}D_+ u_1} \\
& \geq -C \DSN{(u_2)_j} \DSN{u_2} \DSN{u_1} \geq -C \DSNS{u} \DSN{u_j} \\                                                                        & \geq -C \DSNS{u}\left(1+\DSNS{u_j} \right)
\end{align*}
For the third term of (\ref{expan}) we again use the product rule and obtain the expression
\begin{equation*}
\DLIP{3\mu D_+^2\left((f_2)_j  (u_2)_j u_2\right)}{D_+^2 u_1} = \sum_{i_1+i_2+i_3=2}c_{i_1,i_2,i_3} \DLIP{E^{i_2+i_3} D_+^{i_1}(f_2)_j E^{i_3}D_+^{i_2}(u_2)_j D_+^{i_3}u_2}{D_+^2 u_1}
\end{equation*}
For $i_1=2$ and other indices zero we have
\begin{eqnarray*}
\DLIP{(u_2)_j u_2 D_+^2(f_2)_j}{D_+^2 u_1} & \geq & -C\SupNorm{D_+^2 (f_2)_j} \SupNorm{(u_2)_j} \DLN{u_2} \DLN{D_+^2 u_1} \\
                                           & \geq & -C \SupNorm{\frac{d^2}{dx^2}(f_2)_j} \DSN{(u_2)_j} \DSN{u_2} \DSN{u_1} \geq -C \DSN{u_j} \DSNS{u} \\
                                           & \geq & -C \DSNS{u} \left(1+\DSNS{u_j} \right)
\end{eqnarray*}
For $i_2=2$ and other indices zero we have
\begin{eqnarray*}
\DLIP{u_2 E^2 (f_2)_j D_+^2(u_2)_j}{D_+^2 u_1} & \geq & -C\SupNorm{E^2 (f_2)_j} \SupNorm{u_2} \DLN{D_+^2 (u_2)_j} \DLN{D_+^2 u_1} \\
                                               & \geq & -C \DSN{u_2} \DSN{(u_2)_j} \DSN{u_1} \geq -C \DSNS{u} \left( 1+ \DSNS{u_j} \right)
\end{eqnarray*}
For $i_3=2$ and other indices zero we have
\begin{eqnarray*}
\DLIP{E^2 (f_2)_j E^2 (u_2)_j D_+^2 u_2}{D_+^2 u_1} & \geq & -C\SupNorm{E^2 (f_2)_j} \SupNorm{E^2 (u_2)_j} \DLN{D_+^2 u_2} \DLN{D_+^2 u_1} \\
                                                    & \geq & -C \DSN{(u_2)_j} \DSN{u_2} \DSN{u_1} \geq -C \DSNS{u} \left(1+\DSNS{u_j} \right)
\end{eqnarray*}
For $i_1=i_2=1$ and $i_3=0$ we have
\begin{eqnarray*}
\DLIP{u_2 E D_+(f_2)_j D_+ (u_2)_j}{D_+^2 u_1} & \geq & -C\SupNorm{E D_+ (f_2)_j} \SupNorm{u_2} \DLN{D_+ (u_2)_j} \DLN{D_+^2 u_1}\\
                                               & \geq & -C \DSN{u_2} \DSN{(u_2)_j} \DSN{u_1} \geq -C \DSNS{u} \left(1+\DSNS{u_j} \right)
\end{eqnarray*}
For $i_1=i_3=1$ and $i_2=0$ we have
\begin{eqnarray*}
\DLIP{E D_+ (f_2)_j E(u_2)_j D_+ u_2}{D_+^2 u_1} & \geq & -C\SupNorm{E D_+ (f_2)_j} \SupNorm{E (u_2)_j} \DLN{D_+ u_2} \DLN{D_+^2 u_1} \\
                                                 & \geq & -C \DSN{(u_2)_j} \DSN{u_2} \DSN{u_1} \geq -C \DSNS{u} \left(1+\DSNS{u_j} \right)
\end{eqnarray*}
For $i_2=i_3=1$ and $i_1=0$ we have
\begin{eqnarray*}
\DLIP{E^2 (f_2)_j E D_+ (u_2)_j D_+ u_2}{D_+^2 u_1} & \geq & -C\SupNorm{E^2 (f_2)_j} \SupNorm{E D_+ (u_2)_j} \DLN{D_+ u_2} \DLN{D_+^2 u_1} \\
                                                    & \geq & -C \DSN{(u_2)_j} \DSN{u_2} \DSN{u_1} \geq -C \DSNS{u} \left(1+\DSNS{u_j} \right)
\end{eqnarray*}
\indent One can similarly bound all other terms on the left side of (\ref{qbound}) from the expression for $Q_j$ and $\DSIP{\cdot}{\cdot}$ by using the same elementary inequalities as shown above.  
$\square$ 

\begin{lem} \label{lem2.6}
Suppose $K >0$ and $u_0 = \left( (u_1)_0,(u_2)_0\right)$ where $(u_l)_0 \in S^{-\infty}(\mathbb{R} \to \mathbb{R})$ for $l=1,2$ satisfies the property that $\DSN{u_0 } \leq K$ for each $h \in (0,1)$.  Then there exists $T,L, \epsilon >0$ depending only on $K$ such that if $k \in (0, \epsilon)$ and $h\in (0,1)$ then the difference scheme (\ref{dgNLS})  may be solved for each mesh function $u_j = \left( (u_1)_j,(u_2)_j\right)$ with $t_j\in \left[0,T\right]$.  Moreover, we have that $\DSN{ u_j } \leq L$ for each $j$ where $t_j \in \left[0, T\right]$.  \\
\end{lem}

\noindent \textbf{Proof of Lemma 2.6} \hspace{1mm} Choose $T_0 >0$ arbitrarily.  Assume for now that the mesh functions $u_j = \left( (u_1)_{j}, (u_2)_{j} \right)$ are known for each $h,k\in \left(0,1\right)$ and for $0 \leq t_j \leq T_0$.  We will first construct the aformentioned $T,L$ and an $ \epsilon_0 >0$ and show that the mesh functions $u_j$ whose time mesh size satisfies $k \in (0, \epsilon_0) $ will satisfy the inequality $\DSN{ u_j } \leq L$ for each $0\leq t_j \leq T$. \\
\indent Assume that $0\leq t_j, t_{j+1} \leq T$.  Taking inner product of (\ref{dgNLS}) with $u_{j+1}$ we obtain 

\begin{eqnarray*}
 \DSIP{\left(I+kQ_j\right)u_{j+1}}{u_{j+1}} = \DSIP{u_j}{u_{j+1}} - k\mu \DSIP{\left( (f_1)_j(f_2)_j(u_1)_j + (f_1)_j^2(u_2)_j + 3(f_2)_j^2(u_2)_j \hspace{2mm}, \right. \right.\nonumber \\
 \left. \left. -(f_2)_j^2(u_1)_j -3(f_1)_j^2(u_1)_j-2(f_1)_j(f_2)_j(u_2)_j\right) }{u_{j+1}} -k \DSIP{\left( (g_2)_j, -(g_1)_j \right)}{u_{j+1}} 
\end{eqnarray*}

and we may use Cauchy-Schwarz inequality on the right side and simply rewrite the left side to obtain 
\begin{equation}
\DSNS{ u_{j+1} } + k\DSIP{Q_j u_{j+1}}{u_{j+1}} \leq \DSN{ u_{j+1}} \cdot \DSN{ u_j } + k C_{f,g, T_0} \cdot \DSN{u_j}\cdot\DSN{ u_{j+1} }+ k C_{g,T_0}\DSN{ u_{j+1} } \label{L5.1}
\end{equation}
By lemma 2.5 we may choose $C>0$ such that 
\begin{equation}
\DSNS{ u_{j+1}} + k\DSIP{Q_j u_{j+1} }{u_{j+1}} \geq \DSNS{ u_{j+1} } \left[ 1-kC(\DSNS{ u_j } +1) \right] \nonumber
\end{equation}
 By combining (\ref{qbound}) and $( \ref{L5.1} )$ we thus obtain for $t_j,t_{j+1} \in \left[ 0,T_0 \right]$ that 
\begin{eqnarray*}
\DSNS{ u_{j+1} } \big[ 1-kC(\DSNS{ u_j } +1) \big] \leq \DSN{ u_{j+1} } \left( \DSN{ u_j} + k C \DSN{u_j} +k C\right)
\end{eqnarray*}
or equivalently, 
\begin{equation}
\frac{\DSN{ u_{j+1} } - \DSN{ u_j }}{k} \leq C \left( \DSNS{ u_j } + 1\right) \DSN{ u_{j+1} } + C\left( 1+ \DSN{u_j}\right) \nonumber
\end{equation}
 which is an inequality of the form (\ref{dgron}).  Then by lemma 2.4 there exists $0<T\leq T_0$ and $L,\epsilon_0 >0$ depending  on $K \geq \DSN{u_0 }$ such that if $k \in (0, \epsilon_0) $ then $\DSN{ u_j } \leq L$ for each $j$ where $0 \leq t_j\leq T$.  The $T,L,\epsilon_0$ are independent of $h\in (0,1)$ because the constant $C$ is independent of $h$.  Moreover, the $T, L,$ and $\epsilon_0$ depend only on $K, C, P(v):=v^2+1, $ and $Q(v) :=v+1$ by lemma 2.4.  Since $K$ is given and $C$ is determined by the given functions $f$ and $g$ and on the value of $T_0$ we may construct $T,L,$ and $\epsilon_0$ without assuming that $u_j$ is constructed for $0 \leq t_j \leq T_0$.  \\
\indent Given $T,L$ as constructed above it suffices to show that there exists $ 0 < \epsilon \leq \epsilon_0$ such that for any $h\in(0,1)$ and $k \in (0,\epsilon)$ the difference scheme (\ref{dgNLS}) may be solved for $u(x_n,t_j)$ where $ (x_n,t_j)\in \mathbb{R}_h\times \left(\mathbb{R}_k \cap \left[0,T\right] \right)$ - the desired bound $\DSN{ u(\cdot,t_j) } \leq L$ would follow automatically by our construction of $T,L,$ and $\epsilon_0$ given above.\\
\indent Choose $\epsilon >0$ so that $\epsilon C (L^2 +1) < \frac{1}{2}$ and $0<\epsilon \leq \epsilon_0$ and fix values for $h \in (0,1)$ and $k\in (0,\epsilon)$.  
\noindent  Suppose $u_0,u_1,\ldots,u_j$ are known for some $j \geq 0$.  We will show that one may construct $u_{j+1}$ as long as $t_{j+1} \leq T$.  Define an operator $P_j$ := $I+kQ_j$.  Then by lemma 2.5 we have for any mesh function $u=u(x_n)$ 
\begin{equation}
   \DSIP{P_j u}{u}  \geq  \left[1-kC(\DSNS{ u_j }+1) \right]\DSNS{ u } 
     \geq  \left[ 1-\epsilon C(L^2+1)\right]\DSNS{ u } 
     \geq  \frac{1}{2} \DSNS{ u } \label{5.4}
\end{equation}
from which it easily follows that $P_j$ is injective as an operator on mesh functions.  It is also clear that for any mesh function $\rho = \left( \rho_1,\rho_2 \right)\in S_h$ we have $\DSN{P_j \rho} \leq C_{f,g,h,k,u_j, T_0} \DSN{\rho}$ so that $P_j$ is a bounded operator on $S_h$.  In order to solve the difference scheme (\ref{dgNLS}) it is enough to show that the operators $P_j$ are surjective because then $P_j : S_h \rightarrow S_h$ would be a bijection so that by using (\ref{dgNLS}) we could define 
\begin{eqnarray}
u_{j+1} := P_j^{-1} \left[ \left( (u_1)_{j}, (u_2)_{j} \right) - k\mu \left( (f_1)_j(f_2)_j(u_1)_j + (f_1)_j^2(u_2)_j + 3(f_2)_j^2(u_2)_j \hspace{2mm}, \right. \right.  \nonumber \\
  \left. \left. -(f_2)_j^2(u_1)_j -3(f_1)_j^2(u_1)_j-2(f_1)_j(f_2)_j(u_2)_j\right) -k \left( (g_2)_j, -(g_1)_j \right)\right] \label{5.5}
\end{eqnarray}
for $t_{j+1} \leq T$.  First we shall prove that the image of $P_j$ is closed.\\

\noindent \textit{Claim 1} \hspace{1mm} The image of $P_j$: $S_h \rightarrow S_h $ is closed.   \\
\noindent \textit{proof of claim 1} \hspace{1mm} Suppose $P_j u_n \rightarrow v$ as $n\rightarrow \infty$.  We will construct $u\in S_h$ such that $P_j u = v$.  \\
\noindent Inequality (\ref{5.4}) with $u_n-u_m$ implies, by using the Cauchy-Schwarz inequality on the left side, that for any $m,n\in \mathbb{N}$ we have $$ \DSN{P_j (u_n-u_m)} \geq \frac{1}{2} \DSN{u_n-u_m} $$ which implies that the sequence $u_n$ is Cauchy, therefore by completeness of $S_h$ the sequence $u_n$ converges to some $u\in S_h$.  Moreover,
\begin{eqnarray*}
\DSN{P_j u - v} & \leq & \DSN{P_j u - P_j u_n} + \DSN{P_j u_n - v} \leq \left\|P_j \right\| \DSN{u-u_n}+ \DSN{P_j u_n - v}
\end{eqnarray*}  
\noindent and hence for $n$ sufficiently large the right side can be made arbitrarily small so that $P_j u = v$.  This concludes the proof of claim 1.  \\
\noindent \textit{Claim 2} \hspace{1mm}  The operators $P_j $: $S_h \rightarrow S_h $ are surjective.   \\
\noindent \textit{proof of claim 2} \hspace{1mm}  Since the image of $P_j$ is closed we have that $S_h = Im (P_j) \oplus Im (P_j)^{\bot_{S_h}}$.  Suppose there exists $v \in Im (P_j)^{\bot_{S_h}}$.  Then $0 =  \DSIP{P_j v}{v}   \geq \frac{1}{2} \DSNS{v}  $ which implies that $v\equiv 0$.  This proves surjectivity.   \\                             
\noindent Therefore $P_j : S_h \rightarrow S_h$ is a bijection and hence we may solve the difference scheme (\ref{dgNLS}) by defining $u_{j+1}$ as in (\ref{5.5}) as long as $t_j \leq T$ (recall that $(\ref{5.4})$ fails for $t_j > T$ so that $P_j$ would not be invertible after time $T$) however the desired bound $\DSN{u_{j+1}}\leq L$ would be true only if $t_{j+1} \leq T$ . $\square$ 
\begin{rem} 
A similar proof would also work even if the operators $I+kQ_j$ were unbounded as in the case of the KdV and modified KdV equation (see \cite{Bond,Gonz,Menik}).
\end{rem}
\section{Estimates for the Discrete Solutions}
\subsection{Schwartz Boundedness of Discrete Solutions}
\noindent In this section we will show that the solutions to the discrete equation (\ref{dgNLS}) constructed in lemma 2.6 are bounded in all discrete Schwartz norms $\DLN{\IHN{x }^N D^n_+ \cdot}$.  This will follow by some induction arguments shown in next three lemmas.  

\begin{lem}  Let $u_0 = \left( (u_1)_0, (u_2)_0 \right)$ where $(u_l)_0 \in S^{-\infty}(\mathbb{R})$ for $l=1,2$ and let $K,T,\epsilon>0$ all be given as in lemma 2.6 and let $N,n \in \mathbb{N}$.   Suppose that there exists a constant $C_{N,n}>0$ independent of $h \in (0,1)$, $k \in (0,\epsilon)$, and $t_j\in \left[0,T\right]$ such that if $u(x_n,t_j)=\left( u_1 (x_n,t_j),u_2(x_n,t_j) \right)$ is the solution of the difference scheme (\ref{dgNLS}) defined on $\mathbb{R}_h \times \left( \mathbb{R}_k \cap \left[0,T\right] \right)$ with initial condition $u_0 $ then for $0 \leq t_j,t_{j+1} \leq T$ we have 
\begin{eqnarray}
 \DLIP{ \IHN{x}^{2N} D^n_+ (Q_j u_{j+1})}{D^n_+ u_{j+1} } & \geq &  -C_{N,n} \left( \DLNS{ \IHN{x}^N D^n_+ u_{j+1}} + \DLNS{\IHN{x}^N D^n_+ u_j} + 1 \right) \hspace{.5cm} \label{L6.1}
\end{eqnarray}
Then there exists $0 < k_0 \leq \epsilon$ depending on $N$ and $n$ and there exists $C_{N,n}>0$ such that if $h \in (0,1)$, $k\in (0,k_0)$, and $u(x_n,t_j)$ is the solution of the difference scheme (\ref{dgNLS}) defined on $\mathbb{R}_h \times \left( \mathbb{R}_k \cap \left[0,T\right] \right)$ with initial condition $u_0 $ then we have $\DLN{\IHN{x}^N D^n_+u_j}<C_{N,n}$ for $t_j \in \left[0, T\right]$.  Moreover, the constant $C_{N,n}$ is independent of the choice of $h$ and $k$. \\
\end{lem}
\noindent \textbf{Proof of Lemma 3.1} \hspace{1mm} First let us fix values for $h\in (0,1)$ and $k \in (0,\epsilon)$ and let $u(x_n,t_j)$ be the solution of the difference scheme (\ref{dgNLS}) defined on $\mathbb{R}_h \times \left( \mathbb{R}_k \cap \left[0,T\right] \right)$ with initial condition $u_0 $.  If $t_j, t_{j+1}\in \left[0,T\right]$ then we may apply $D_+^n$ to both sides of (\ref{dgNLS}) and take the inner product $\DLIP{\IHN{x}^{2N}\cdot}{\cdot}$ of the resulting equation with $D_+^n u_{j+1}$ and use proposition 2.1, the product rule, and corollary 2.3 to obtain 

\begin{align}
& \DLIP{\IHN{x}^{2N}D_+^n \big(I+kQ_j \big) u_{j+1}}{D_+^n u_{j+1}}  =   \DLIP{\IHN{x}^{2N}D_+^n u_j}{D_+^n u_{j+1}} \hspace{6cm}\nonumber \\
                                                                   &    - k\mu \DLIP{\IHN{x}^{2N}D_+^n \left( (f_1)_j(f_2)_j(u_1)_j + (f_1)_j^2(u_2)_j + 3(f_2)_j^2(u_2)_j , \right. \right. \nonumber \\
                                                                   &   \left. \left. -(f_2)_j^2(u_1)_j -3(f_1)_j^2(u_1)_j-2(f_1)_j(f_2)_j(u_2)_j\right) }{D_+^n u_{j+1}}
                                                                      - k \DLIP{\IHN{x}^{2N}D_+^n \left( (g_2)_j, -(g_1)_j \right)}{D_+^n u_{j+1}}  \nonumber \\
                                                                   & \leq  \DLN{\IHN{x}^N D_+^n u_{j+1}} \DLN{\IHN{x}^N D_+^n u_j} +kC_{f,g,K}\DLN{\IHN{x}^N D_+^n u_{j+1}}\left(1+ \DLN{\IHN{x}^N D_+^n u_j} \right) \nonumber \\
                                                                   & \leq  \frac{1}{2}\DLNS{\IHN{x}^N D_+^n u_{j+1}} + \frac{1}{2} \DLNS{\IHN{x}^N D_+^n u_j} + k C_{N,n,f,g}  \left( 1+ \DLNS{\IHN{x}^N D_+^n u_{j}}+ \DLNS{\IHN{x}^N D_+^n u_{j+1}} \right)  \label{L6.2}
\end{align}

\noindent We may then continue the left side of $(\ref{L6.2})$ by using $(\ref{L6.1})$ to obtain 
\begin{align}
& \DLIP{\IHN{x}^{2N}D_+^n \big(I+kQ_j \big) u_{j+1}}{D_+^n u_{j+1}} = \DLNS{\IHN{x}^N D_+^n u_{j+1}} + k\DLIP{\IHN{x}^{2N}D_+^n Q_j u_{j+1}}{D_+^n u_{j+1}} \hspace{22mm} \nonumber \\  
 & \geq \DLNS{\IHN{x}^N D_+^n u_{j+1}} - k C \left( \DLNS{\IHN{x}^N D_+^n u_j} + \DLNS{\IHN{x}^N D_+^n u_{j+1}} + 1 \right) \label{L6.3}
\end{align}
By combining $(\ref{L6.2})$ and $(\ref{L6.3})$ we obtain
\begin{eqnarray*}
 \frac{1}{k}\left(\DLNS{\IHN{x}^N D_+^n u_{j+1}} - \DLNS{\IHN{x}^N D_+^n u_j} \right) 
 \leq  C \left( \DLNS{\IHN{x}^N D_+^n u_{j+1}} + \DLNS{\IHN{x}^N D_+^n u_j}  + 1 \right)
\end{eqnarray*}
which is an inequality of the form (\ref{dgron}) where the right side is linear.  Now we may invoke lemma 2.4 to obtain the existence of $k_0 >0$ and $C_{N,n}>0$ such that if $k \in (0,k_0)$ and $t_j \in \left[0,T\right]$ then $$\DLN{\IHN{x}^N D_+^n u_j} < C_{N,n} $$ and since $C$ is independent of $h$ and $k$ we also have that $C$ is independent of $h$ and $k$.  Thus we have proven that for each $0<h\leq 1$, $0<k\leq k_0$, $0 \leq t_j \leq T$ and solution $u(x_n,t_j)$ of (\ref{dgNLS}) defined on $\mathbb{R}_h\times\left( \mathbb{R}_k \cap \left[0,T\right]\right)$ with initial data $u_0$ we have the inequality $\DLN{\IHN{x}^N D_+^n u_j} < C_{N,n} $.  Moreover since $C$ depends on $N$ and $n$ it follows that $k_0$ also depends on $N$ and $n$. $\square$ 

\begin{lem}
Let $u_0= \left( (u_1)_0, (u_2)_0 \right) $ where $(u_l)_0 \in S^{-\infty}(\mathbb{R})$ for $l=1,2$ and let $K,T,L,\epsilon>0$ all be given as in lemma 2.6 and let $n \in \mathbb{N}$.  Then there exists $0 < k_0 \leq \epsilon$ depending on $n$ and there exists $C_n >0$ such that if $h \in (0,1)$, $k\in (0,k_0)$, and $u(x_n,t_j) = \left( u_1 (x_n,t_j),u_2(x_n,t_j) \right)$ is the solution of the difference scheme (\ref{dgNLS}) defined on $\mathbb{R}_h \times \left( \mathbb{R}_k \cap \left[0,T\right] \right)$ with initial condition $u_0 $ then we have $\DLN{\IHN{x} D^n_+u_j}<C_n$ for $t_j \in \left[0, T\right]$.  Moreover, the constant $C_n$ is independent of the choice of $h$ and $k$. 
\end{lem}

\noindent \textbf{Proof of Lemma 3.2} \hspace{1mm} By lemma 2.6 the statement is true for $n = 0, 1$ by taking $k_0 = \epsilon$ and $C_n = L$.  We shall prove the statement by induction on $n$.  Assume it is true for all $l\leq n-1$.  We will prove that it is true for $l=n$.  First we shall use the inductive hypothesis to prove some slightly weaker claims which we give below. \\

\noindent \textit{Claim 1} \hspace{1mm} There exists $0<\epsilon_0 \leq \epsilon$ and $C_n>0$ such that if $h\in (0,1)$ and $k\in (0,\epsilon_0)$ and $u(x_n,t_j)$ is the solution of the difference scheme (\ref{dgNLS}) defined on $\mathbb{R}_h \times \left( \mathbb{R}_k \cap \left[0,T\right] \right)$ with initial condition $u_0 $ then we have $\DLN{D_+^n u_j} \leq C_n$  for $t_j \in \left[0,T\right]$. \\
\noindent \textit{proof of claim 1} \hspace{1mm} By lemma 3.1 it suffices to prove that there exists a constant $C_n>0$ independent of $h \in (0,1)$, $k \in (0,\epsilon)$, and $t_j\in \left[0,T\right]$ such that if $u(x_n,t_j)$ is the solution of the difference scheme (\ref{dgNLS}) defined on $\mathbb{R}_h \times \left( \mathbb{R}_k \cap \left[0,T\right] \right)$ with initial condition $u_0 $ then for $0 \leq t_j,t_{j+1} \leq T$ we have 
\begin{eqnarray*}
 \DLIP{ D^n_+ (Q_j u_{j+1})}{D^n_+ u_{j+1} }  \geq  -C_n \left( \DLNS{ D^n_+ u_{j+1}} + \DLNS{ D^n_+ u_j} + 1 \right) 
\end{eqnarray*}

\noindent To this end we shall fix values for $h\in (0,1)$ and $k\in (0,\epsilon)$ and a solution $u(x_n,t_j)$ of the difference scheme (\ref{dgNLS}) defined on $\mathbb{R}_h \times \left( \mathbb{R}_k \cap \left[0,T\right] \right)$ with initial condition $u_0 $.  In order to prove that the above estimate holds for some constant $C_n>0$ we simply prove that the estimate can be made for each term of $Q_j$ and then by adding all these estimates we will obtain the estimate for $Q_j$.  These estimates are analogous to those given in lemma 2.5 however we will also use the inductive hypothesis.  For simplicity we will ignore all occurences of the shift $E$.  Here are the necessary estimates. \\

\noindent \underline{\textit{Estimate for Term}} $\DLIP{D_+^n (D_+ D_- (u_2)_{j+1}) }{D_+^n (u_1)_{j+1}} - \DLIP{D_+^n (D_+ D_- (u_1)_{j+1}) }{D_+^n (u_2)_{j+1}}$: 
\begin{eqnarray*}
\DLIP{D_+^n (D_+ D_- (u_2)_{j+1}) }{D_+^n (u_1)_{j+1}} - \DLIP{D_+^n (D_+ D_- (u_1)_{j+1}) }{D_+^n (u_2)_{j+1}}  \\ = -\DLIP{D_+^{n+1} (u_2)_{j+1}}{D_+^{n+1} (u_1)_{j+1}} + \DLIP{D_+^{n+1} (u_1)_{j+1}}{D_+^{n+1} (u_2)_{j+1}} = 0
\end{eqnarray*}

\noindent \underline{\textit{Estimate for Term}} $\DLIP{D_+^n ((u_1)_j(u_2)_j(u_1)_{j+1}) }{D_+^n (u_1)_{j+1}}$ : \\ 
\indent By the product rule we obtain,
\begin{align*}
& \DLIP{D_+^n ((u_1)_j(u_2)_j(u_1)_{j+1}) }{D_+^n (u_1)_{j+1}} \\ 
& = \sum_{i_1 + i_2 + i_3=n} c_{i_1,i_2,i_3} \DLIP{ D_+^{i_1}(u_1)_j D_+^{i_2}(u_2)_j  D_+^{i_3}(u_1)_{j+1} }{D_+^n (u_1)_{j+1}}
\end{align*}
If $ 2 \leq i_1 + i_2 \leq n$  and $i_1\leq i_2$ then we have the estimate,
\begin{align*}
& \DLIP{  D_+^{i_1}(u_1)_j  D_+^{i_2}(u_2)_j  D_+^{i_3}(u_1)_{j+1} }{D_+^n (u_1)_{j+1}} \\
& \geq -C \SupNorm{D_+^{i_1}(u_1)_j} \DLN{D_+^{i_2}(u_2)_j} \SupNorm{D_+^{i_3}(u_1)_{j+1}} \DLN{D_+^n(u_1)_{j+1}} \\
                                                                     & \geq -C \left( \DLN{(u_2)_j} + \DLN{D_+^n(u_2)_j} \right) \DLN{D_+^n(u_1)_{j+1}} \\
                                                                     & \geq -C \left(1+ \DLNS{D_+^n (u_1)_{j+1}} + \DLNS{D_+^n (u_2)_{j}} \right) \\
                                                                     & \geq -C \left(1+ \DLNS{D_+^n u_{j+1}} + \DLNS{D_+^n u_{j}} \right)  
\end{align*}
Also, if instead we had $ 2 \leq i_1 + i_2 \leq n$  and $i_2\leq i_1$ then the same estimates work by just switching $(u_1)_j$ and $(u_2)_j$. 

For the term where $ i_1 + i_2 = 1$  and $i_1=1$ we have the estimate,
\begin{align*}
& \DLIP{ (u_2)_j D_+ (u_1)_j  D_+^{n-1}(u_1)_{j+1} }{D_+^n (u_1)_{j+1}} \\
& \geq -\SupNorm{(u_2)_j} \SupNorm{D_+^{n-1}(u_1)_{j+1}} \DLN{D_+(u_1)_{j}} \DLN{D_+^n(u_1)_{j+1}} \\
                                                                      & \geq -C\left( \DLN{(u_1)_{j+1}} + \DLN{D_+^n (u_1)_{j+1}} \right) \DLN{D_+^n (u_1)_{j+1}} \\
                                                                      & \geq -C \left(1+ \DLNS{D_+^n (u_1)_{j+1}} \right) \geq -C \left(1+ \DLNS{D_+^n u_{j+1}} \right)
\end{align*}
Also, if instead we had $i_1 + i_2 = 1$  and $i_2=1$ then the same estimates work by just switching $(u_1)_j$ and $(u_2)_j$. 

For the term where $i_1 + i_2 = 0$ we have the estimate,
\begin{eqnarray*}
\DLIP{ (u_1)_j (u_2)_j D_+^{n}(u_1)_{j+1} }{D_+^n (u_1)_{j+1}} & \geq & -\SupNorm{(u_1)_j} \SupNorm{(u_2)_j} \DLNS{D_+^n (u_1)_{j+1}} \\
                                                               & \geq & -C \left( \DLN{u_j } + \DLN{D_+ u_j } \right)^2 \DLNS{D_+^n u_{j+1}}  \\
                                                               & \geq & -C \DLNS{D_+^n u_j}
\end{eqnarray*}
\noindent All other terms can be similarly estimated.  This concludes the proof of claim 1.  \\

\noindent \textit{Claim 2} \hspace{1mm}   There exists $0<\epsilon_1 \leq \epsilon_0$ and $C_n>0$ such that if $h\in (0,1)$ and $k\in (0,\epsilon_1)$ and $u(x_n,t_j)$ is the solution of the difference scheme (\ref{dgNLS}) defined on $\mathbb{R}_h \times \left( \mathbb{R}_k \cap \left[0,T\right] \right)$ with initial condition $u_0 $ then we have $\DLN{D_+^{n+1} u_j} \leq C_n$ for $t_j \in \left[0,T\right]$. \\
\noindent \textit{proof of claim 2} \hspace{1mm} As in the proof of claim 1 we see that by lemma 3.1 it suffices to prove that there exists a constant $C_n>0$ independent of $h \in (0,1)$, $k \in (0, \epsilon_0 )$, and $t_j\in \left[0,T\right]$ such that if $u(x_n,t_j)$ is the solution of the difference scheme (\ref{dgNLS}) defined on $\mathbb{R}_h \times \left( \mathbb{R}_k \cap \left[0,T\right] \right)$ with initial condition $u_0 $ then for $0 \leq t_j,t_{j+1} \leq T$ we have 
\begin{eqnarray*}
 \DLIP{ D^{n+1}_+ (Q_j u_{j+1})}{D^{n+1}_+ u_{j+1} }  \geq  -C_n \left( \DLNS{ D^{n+1}_+ u_{j+1}} + \DLNS{ D^{n+1}_+ u_j} + 1 \right) 
\end{eqnarray*}

\noindent To this end we shall again fix values for $h\in (0,1)$ and $k\in (0,\epsilon_0)$ and a solution $u(x_n,t_j)$ to the solution of the difference scheme (\ref{dgNLS}) defined on $\mathbb{R}_h \times \left( \mathbb{R}_k \cap \left[0,T\right] \right)$ with initial condition $u_0 $.  In order to prove that the above estimate holds for some constant $C_n>0$ we simply prove that the estimate can be made for each term of $Q_j$ and then by adding all these estimates we will obtain the estimate for $Q_j$.  By using the inductive hypothesis and the result of claim 1, it is easily seen that the estimates shown in the proof of claim 1 can all be applied with $n$ replaced by $n+1$ and the others can be similarly estimated.  This concludes the proof of claim 2.\\

\indent  Now we return to the proof of the lemma in the case $l=n$.  By using the same reasoning as in
the above claims we see by lemma 3.1 that in order to construct $k_0 \in (0, \epsilon_2)$ and $C_n$ it suffices to
prove the estimate
\begin{eqnarray*}
\DLIP{ \IHN{x}^{2} D^n_+ (Q_j u_{j+1})}{D^n_+ u_{j+1} } & \geq &  -C_{n} \left( \DLNS{ \IHN{x} D^n_+ u_{j+1}} + \DLNS{\IHN{x} D^n_+ u_j} + 1 \right)
\end{eqnarray*}
This estimate will follow by adding all of the below estimates. In the below estimates we will use our inductive hypothesis and the estimates $
\DLN{D_+^{n+i} u_j } \leq C_{n+i}$ for $i=0,1$.  For conciseness we shall again only formulate estimates for a few terms because the rest can be bounded similarly. 

\noindent \underline{\textit{Estimate for Term}} $\DLIP{D_+ D_- D_+^n (u_2)_{j+1} }{\IHN{x}^2 D_+^n (u_1)_{j+1}} - \DLIP{D_+ D_- D_+^n (u_1)_{j+1} }{\IHN{x}^2 D_+^n (u_2)_{j+1}}$:  
\begin{eqnarray*}
\DLIP{D_+ D_- (D_+^n (u_2)_{j+1}) }{\IHN{x}^2 D_+^n (u_1)_{j+1}} - \DLIP{D_+ D_- (D_+^n (u_1)_{j+1}) }{\IHN{x}^2 D_+^n (u_2)_{j+1}} \hspace{1cm}\\
 =  - \DLIP{D_+^{n+1} (u_2)_{j+1}}{E\IHN{x}^2 D_+^{n+1}(u_1)_{j+1}} - \DLIP{D_+^{n+1} (u_2)_{j+1}}{D_+^n (u_1)_{j+1} D_+ \IHN{x}^2} \\
   + \DLIP{D_+^{n+1} (u_1)_{j+1}}{E\IHN{x}^2 D_+^{n+1} (u_2)_{j+1}} + \DLIP{D_+^{n+1} (u_1)_{j+1}}{D_+^n (u_2)_{j+1} D_+ \IHN{x}^2} \\ 
 =   \DLIP{D_+^{n+1} (u_1)_{j+1}}{D_+^n (u_2)_{j+1} D_+ \IHN{x}^2} - \DLIP{D_+^{n+1} (u_2)_{j+1}}{D_+^n (u_1)_{j+1} D_+ \IHN{x}^2} \hspace{5mm} \\
 \geq  - \DLN{D_+^{n+1} (u_1)_{j+1}} \DLN{\IHN{x} D_+^n (u_2)_{j+1}} - \DLN{D_+^{n+1} (u_2)_{j+1}} \DLN{\IHN{x} D_+^n (u_1)_{j+1}} \hspace{2mm} \\
 \geq  -C \left( \DLN{\IHN{x}D_+^n (u_2)_{j+1}} + \DLN{\IHN{x}D_+^n (u_1)_{j+1}} \right)  \geq  -C \left( \DLNS{\IHN{x} D_+^n u_{j+1}} + 1 \right) \hspace{5mm}
\end{eqnarray*}

\noindent \underline{\textit{Estimate for Term}} $\DLIP{\IHN{x} D_+^n ((u_1)_j(u_2)_j(u_1)_{j+1}) }{\IHN{x}D_+^n (u_1)_{j+1}}$ : \\ 
\indent By the product rule we obtain,
\begin{align*}
& \DLIP{\IHN{x} D_+^n ((u_1)_j(u_2)_j(u_1)_{j+1}) }{\IHN{x}D_+^n (u_1)_{j+1}} \\
& = \sum_{i_1 + i_2 + i_3=n} c_{i_1,i_2,i_3} \DLIP{ \IHN{x} D_+^{i_1}(u_1)_j D_+^{i_2}(u_2)_j  D_+^{i_3}(u_1)_{j+1} }{\IHN{x} D_+^n (u_1)_{j+1}}
\end{align*}
If $ 2 \leq i_1 + i_2 \leq n$  and $i_1\leq i_2$ then we have the estimate,
\begin{eqnarray*}
\DLIP{ \IHN{x} D_+^{i_1}(u_1)_j  D_+^{i_2}(u_2)_j  D_+^{i_3}(u_1)_{j+1} }{\IHN{x} D_+^n (u_1)_{j+1}} \hspace{5cm} \\
 \geq  -C \SupNorm{D_+^{i_1}(u_1)_j} \SupNorm{D_+^{i_3} (u_1)_{j+1}} \DLN{\IHN{x}D_+^{i_2} (u_2)_j} \DLN{\IHN{x}D_+^n (u_1)_{j+1}} \\
 \geq  -C \DLN{\IHN{x} D_+^n (u_1)_{j+1}} \geq -C \left(1+ \DLNS{\IHN{x} D_+^n u_{j+1}} \right) \hspace{2.6cm}
\end{eqnarray*}
Also, if instead we had $ 2 \leq i_1 + i_2 \leq n$  and $i_2\leq i_1$ then the same estimates work by just switching $(u_1)_j$ and $(u_2)_j$. 

For the term where $ i_1 + i_2 = 1$  and $i_1=1$ we have the estimate,
\begin{eqnarray*}
\DLIP{\IHN{x} (u_2)_j D_+ (u_1)_j  D_+^{n-1}(u_1)_{j+1} }{\IHN{x} D_+^n (u_1)_{j+1}} \hspace{5cm}\\ 
\geq  - \SupNorm{\IHN{x}(u_2)_j} \SupNorm{D_+ (u_1)_{j}} \DLN{\IHN{x} D_+^{n-1}(u_1)_{j+1}} \DLN{\IHN{x} D_+^n (u_1)_{j+1}} \\
                                                                       \geq  -C \DLN{\IHN{x} D_+^n (u_1)_{j+1}} \geq  -C \left(1+ \DLNS{\IHN{x}D_+^n u_{j+1}} \right) \hspace{2.5cm}
\end{eqnarray*}
Also, if instead we had $i_1 + i_2 = 1$  and $i_2=1$ then the same estimates work by just switching $(u_1)_j$ and $(u_2)_j$. 

For the term where $i_1 + i_2 = 0$ we have the estimate,
\begin{eqnarray*}
\DLIP{ \IHN{x} (u_1)_j (u_2)_j D_+^{n} (u_1)_{j+1} }{\IHN{x} D_+^n (u_1)_{j+1} }  \geq  -\SupNorm{\IHN{x}(u_1)_j} \SupNorm{\IHN{x}(u_2)_j} \DLNS{D_+^n (u_1)_{j+1}} \geq -C
\end{eqnarray*}
\noindent  This concludes the proof of the lemma. \hspace{2cm} $\square$ 

\begin{lem}
Let $u_0 = \left( (u_1)_0, (u_2)_0 \right)$ where $(u_l)_0 \in S^{-\infty}(\mathbb{R})$ for $l=1,2$ and let $K,T,\epsilon>0$ all be given as in lemma 2.6 and let $N,n \in \mathbb{N}$.  Then there exists $0 < k_0 \leq \epsilon$ depending on $N$ and $n$ and there exists $C_{N,n} >0$ such that if $h \in (0,1)$, $k\in (0,k_0)$, and $u(x_n,t_j) = \left( u_1 (x_n,t_j),u_2(x_n,t_j) \right)$ is the solution of the difference scheme (\ref{dgNLS}) defined on $\mathbb{R}_h \times \left( \mathbb{R}_k \cap \left[0,T\right] \right)$ with initial condition $u_0 $ then we have $\DLN{\IHN{x}^N D^n_+u_j}<C_{N,n}$ for $t_j \in \left[0, T\right]$.  Moreover, the constant $C_{N,n}$ is independent of the choice of $h$ and $k$. 
\end{lem}

\noindent \textbf{Proof of Lemma 3.3} \hspace{1mm} By lemma 3.2 the assertion is true for $N=1$ and for all $n\in\mathbb{N}$.  We prove the assertion by induction on $N$.  Assume it is true for all $M\leq N-1$ and for all $n\in\mathbb{N}$.  We will show that for $M=N$ the statement is satisfied by all $n\in \mathbb{N}$ by induction on $n$.  \\
\indent  Let us denote the value of $k_0$ corresponding to a particular value of $N$ and $n$ by $k_0(N,n)$.  By construction of $k_0(N,n)$ for its known values we see that if $N_2 < N_1$ then $k_0(N_1,n)\leq k_0(N_2,n)$ and if $n_2 < n_1$ then $k_0(N,n_1)\leq k_0(N,n_2)$.  Let $\epsilon_0 := k_0(N-1,1)$.  By use of lemma 3.1 and taking $\epsilon = \epsilon_0$ our statement for $n=0$ will follow if we prove that there exists a constant $C_{N,0}>0$ independent of $h \in (0,1)$, $k \in (0,\epsilon_0)$, and $t_j\in \left[0,T\right]$ such that if $u(x_n,t_j)$ is the solution of the difference scheme (\ref{dgNLS}) defined on $\mathbb{R}_h \times \left( \mathbb{R}_k \cap \left[0,T\right] \right)$ with initial condition $u_0 $ then for $0 \leq t_j,t_{j+1} \leq T$ we have

\begin{eqnarray*}
\DLIP{ \IHN{x}^{2N}  Q_j u_{j+1}}{ u_{j+1} } & \geq &  -C_{N,0} \left( \DLNS{ \IHN{x}^N u_{j+1}} + \DLNS{\IHN{x}^N u_j} + 1 \right)
\end{eqnarray*}

\noindent To this end we shall fix values for $h\in (0,1)$ and $k\in (0,\epsilon_0)$ and a solution $u(x_n,t_j)$ of the difference scheme (\ref{dgNLS}) defined on $\mathbb{R}_h \times \left( \mathbb{R}_k \cap \left[0,T\right] \right)$ with initial condition $u_0 $.  In order to prove that the above estimate holds for some constant $C_{N,0}>0$ we simply prove that the estimate can be made for each term of $Q_j$ and then by adding all these estimates we will obtain the estimate for $Q_j$.  We will again only show a few estimates since all can be done similarly and we will ignore all occurences of the shift $E$ for simplicity.\\

\noindent \underline{\textit{Estimate for Term}} $\DLIP{D_+ D_- (u_2)_{j+1} }{\IHN{x}^{2N}(u_1)_{j+1}} - \DLIP{D_+ D_- (u_1)_{j+1} }{ \IHN{x}^{2N}(u_2)_{j+1}}$:  
\begin{align*}
& \DLIP{D_+ D_- (u_2)_{j+1} }{\IHN{x}^{2N}(u_1)_{j+1}} - \DLIP{D_+ D_- (u_1)_{j+1} }{ \IHN{x}^{2N}(u_2)_{j+1}}\\  
& =  \DLIP{D_+ (u_1)_{j+1} }{(u_2)_{j+1} D_+ \IHN{x}^{2N}} - \DLIP{D_+ (u_2)_{j+1} }{(u_1)_{j+1} D_+ \IHN{x}^{2N}}  \\  
& \geq  -C \DLN{\IHN{x}^{N-1}D_+ (u_1)_{j+1}} \DLN{\IHN{x}^N (u_2)_{j+1}} -C \DLN{\IHN{x}^{N-1} D_+ (u_2)_{j+1}} \DLN{\IHN{x}^N (u_1)_{j+1}}\\ 
& \geq -C \DLN{\IHN{x}^N u_{j+1}}  
\end{align*}

\noindent \underline{\textit{Estimate for Term}} $\DLIP{\IHN{x}^N (u_1)_j(u_2)_j(u_1)_{j+1} }{\IHN{x}^N (u_1)_{j+1}}$:
\begin{eqnarray*}
\DLIP{\IHN{x}^N (u_1)_j(u_2)_j(u_1)_{j+1} }{\IHN{x}^N (u_1)_{j+1}} & \geq & \SupNorm{(u_1)_j} \SupNorm{(u_2)_j} \DLNS{\IHN{x}^N (u_1)_{j+1}}\\
                                                                   & \geq &-C \DLNS{\IHN{x}^N u_{j+1}}
\end{eqnarray*}

\noindent This concludes the proof for the case $n=0$.  \\

\indent Now we will assume that the statement is true for all $l\leq n-1$ and we will prove that it is true for $l=n$.  Let $\epsilon_1 := k_0(N-1,n+1)$.  The proof will again follow from lemma 3.1 by taking $\epsilon = \epsilon_1$ if we can prove that there exists a constant $C_{N,n}>0$ independent of $h \in (0,1)$, $k \in (0,\epsilon_1)$, and $t_j\in \left[0,T\right]$ such that if $u(x_n,t_j)$ is the solution of the difference scheme (\ref{dgNLS}) defined on $\mathbb{R}_h \times \left( \mathbb{R}_k \cap \left[0,T\right] \right)$ with initial condition $u_0 $ then for $0 \leq t_j,t_{j+1} \leq T$ we have
$$
 \DLIP{ \IHN{x}^{2N} D^n_+ (Q_j u_{j+1})}{D^n_+ u_{j+1} }   \geq    -C_{N,n} \left( \DLNS{ \IHN{x}^N D^n_+ u_{j+1}} + \DLNS{\IHN{x}^N D^n_+ u_j} + 1 \right) 
$$

\noindent  To this end we shall fix values for $h\in (0,1)$ and $k\in (0,\epsilon_1)$ and a solution $u(x_n,t_j)$ of the difference scheme (\ref{dgNLS}) defined on $\mathbb{R}_h \times \left( \mathbb{R}_k \cap \left[0,T\right] \right)$ with initial condition $u_0 $.  In order to prove that the above estimate holds for some constant $C_{N,n}>0$ we again prove that the estimate can be made for each term of $Q_j$ and then by adding all these estimates we will obtain the estimate for $Q_j$.  Here are estimates for some terms - the rest can be estimated similarly. \\

\noindent \underline{\textit{Estimate for Term}} $\DLIP{D_+ D_- D_+^n (u_2)_{j+1} }{\IHN{x}^{2N}D_+^n (u_1)_{j+1}} - \DLIP{D_+ D_- D_+^n (u_1)_{j+1} }{ \IHN{x}^{2N} D_+^n (u_2)_{j+1}}$:  
\begin{align*}
& \DLIP{D_+ D_- D_+^n (u_2)_{j+1} }{\IHN{x}^{2N}D_+^n (u_1)_{j+1}} - \DLIP{D_+ D_- D_+^n (u_1)_{j+1} }{ \IHN{x}^{2N}D_+^n (u_2)_{j+1}}  \\ 
& =  \DLIP{D_+^{n+1} (u_1)_{j+1} }{D_+^n (u_2)_{j+1} D_+ \IHN{x}^{2N}} - \DLIP{D_+^{n+1} (u_2)_{j+1} }{D_+^n (u_1)_{j+1} D_+ \IHN{x}^{2N}} \\  
& \geq  -C \DLN{\IHN{x}^{N-1}D_+^{n+1} (u_1)_{j+1}} \DLN{\IHN{x}^N D_+^n (u_2)_{j+1}} \\
& -C \DLN{\IHN{x}^{N-1} D_+^{n+1} (u_2)_{j+1}} \DLN{\IHN{x}^N D_+^n (u_1)_{j+1}}\\ 
& \geq -C \left( \DLN{\IHN{x}^N D_+^n (u_1)_{j+1}} + \DLN{\IHN{x}^N D_+^n (u_2)_{j+1}} \right) \\
& \geq -C \left(1+\DLNS{\IHN{x}^N D_+^n u_{j+1}} \right)  
\end{align*}

\noindent \underline{\textit{Estimate for Term}} $\DLIP{\IHN{x}^{2N} D_+^n ((u_1)_j(u_2)_j(u_1)_{j+1}) }{ D_+^n (u_1)_{j+1}}$ : \\ 
\indent For simplicity we will ignore the shifts E.  By the product rule we obtain,
\begin{align*}
& \DLIP{\IHN{x}^{2N} D_+^n ((u_1)_j(u_2)_j(u_1)_{j+1}) }{ D_+^n (u_1)_{j+1}} \\
& = \sum_{i_1 + i_2 + i_3=n} c_{i_1,i_2,i_3} \DLIP{ \IHN{x}^{2N} D_+^{i_1}(u_1)_j D_+^{i_2}(u_2)_j  D_+^{i_3}(u_1)_{j+1} }{ D_+^n (u_1)_{j+1}}
\end{align*}
For each term we have the estimate,
\begin{align}
& \DLIP{ \IHN{x}^N D_+^{i_1}(u_1)_j  D_+^{i_2}(u_2)_j  D_+^{i_3}(u_1)_{j+1} }{\IHN{x}^N D_+^n (u_1)_{j+1}}\nonumber \\
& \geq  -C \SupNorm{D_+^{i_1}(u_1)_j} \SupNorm{D_+^{i_2} (u_2)_{j}} \DLN{\IHN{x}^N D_+^{i_3} (u_1)_{j+1}} \DLN{\IHN{x}^N D_+^n (u_1)_{j+1}} \nonumber \\
& \geq  -C \left( \DLN{\IHN{x}^N (u_1)_{j+1}} + \DLN{\IHN{x}^N D_+^n (u_1)_{j+1}} \right) \DLN{\IHN{x}^N D_+^n (u_1)_{j+1}} \nonumber \\
& \geq -C \left( 1+ \DLNS{\IHN{x}^N D_+^n u_{j+1}} \right) \nonumber
\end{align}

\noindent This concludes the proof of the lemma. $\square$ 

\begin{cor}
Let $u_0 = \left( (u_1)_0, (u_2)_0 \right)$ where $(u_l)_0 \in S^{-\infty}(\mathbb{R})$ for $l=1,2$ and let $K,T,\epsilon>0$ all be given as in lemma 2.6 and let $N,n \in \mathbb{N}$.  Then there exists $0 < k_0 \leq \epsilon$ depending on $N$ and $n$ and there exists $C_{N,n} > 0$ such that if $h \in (0,1)$, $k\in (0,k_0)$, and $u(x_n,t_j)= \left( u_1 (x_n,t_j),u_2(x_n,t_j) \right)$ is the solution of the difference scheme (\ref{dgNLS}) defined on $\mathbb{R}_h \times \left( \mathbb{R}_k \cap \left[0,T\right] \right)$ with initial condition $u_0 $ then we have $\DLN{D^n_+ \left(x^N u_j\right)}<C_{N,n}$ for $t_j \in \left[0, T\right]$.  Moreover, the constant $C_{N,n}$ is independent of the choice of $h$ and $k$. 
\end{cor}
\noindent \textbf{Proof of Corollary 3.4} \hspace{1mm} By the product rule for $D_+$ this statement can be proven by induction on $n$.  We shall omit the necessary details here. $\square$
\subsection{Boundedness of Time-Differentiated Extended Discrete Solutions} 
\noindent In this section we will show that a certain time-extension of the discrete solution with domain $\mathbb{R}_h \times \mathbb{R}_k$ remains bounded in the Schwartz semi-norms $\DLN{D_{t,+}^m D_+^n (x^N \cdot)}$ for $n,N \in \mathbb{N}$ and $m\in \left\{0,1,2,3\right\}$. \\
\indent Suppose $u(x_n,t_j)=\left(u_1(x_n,t_j),u_2(x_n,t_j)\right)$ is mesh function defined on $\mathbb{R}_h \times \left( \mathbb{R}_k \cap \left[0,T\right] \right)$ for some $h,k,T>0$ where $h\in (0,1)$ and $k \leq T/3$.  We shall define an extension of $u$ to $\mathbb{R}_h \times \mathbb{R}_k$ by the following: \\
Let $\phi(t) \in C^{\infty}_c\left(\mathbb{R}\right)$ such that $\phi(t) = 1$ for $t \in \left[-1,T+1\right]$ and $\phi(t) = 0$ for $t \notin \left[-2,T+2\right]$.  \\
For $t_j > T$ we define recursively 

$$u(x_n,t_j) :=  u(x_n,t_{j-1}) + k D_{t,+} u(x_n,t_{j-2}) + k^2 D_{t,+}^2 u(x_n,t_{j-3}) $$ 

\noindent and similarly for $t_j <0$ we define 

$$u(x_n,t_j) :=  u(x_n,t_{j+1}) - k D_{t,+} u(x_n,t_{j+1})+k^2 D_{t,+}^2 u(x_n,t_{j+1})  $$ 

\noindent Then we define $\hat{u}:= \phi \cdot u $, which is clearly an extension of $u$ to $\mathbb{R}_h \times \mathbb{R}_k$ and is compactly supported in time.  \\
\indent It is possible to prove that a finite time discrete solution extended in this way remains bounded in the discrete Schwartz semi-norms $\DLN{D_{t,+}^m D^n_+ \left( x^N u \right)}$ for $m=0,1,2,3$ (see lemma 3.5 below).  The proof is omitted here but the same statement appears in \cite{Gonz} and is proved there in lemma 3.5.  It relies primarily on the estimates of the type given in corollary 3.4 above.  
\begin{lem}
Let $u_0=\left( (u_1)_0,(u_2)_0 \right)$ where $(u_l)_0 \in S^{-\infty}(\mathbb{R})$ for $l=1,2$ and let $K,T,\epsilon >0$ all be given as in lemma 2.6 and let $N,n \in \mathbb{N}$, $m \in \left\{0,1,2,3\right\}$.  Then there exists $0 < k_0 \leq \epsilon$ depending on $m,n$ and $N$ and there exists $C_{m,N,n} > 0$ such that if $h \in (0,1)$, $k\in (0,k_0)$, and $u(x_n,t_j)$ is the solution of the difference scheme (\ref{dgNLS}) defined on $\mathbb{R}_h \times \left( \mathbb{R}_k \cap \left[0,T\right] \right)$ with initial condition $u_0 $ then we have 
\begin{equation}
\DLN{D_{t,+}^m D^n_+ \left( x^N \hat{u}(x, t_j) \right)}\leq C  \label{L9.1}
\end{equation}  
\noindent for each $t_j \in \mathbb{R}_k$.
\end{lem}

\section{Obtaining Smooth Solutions from Discrete Solutions}
\subsection{The Smoothing Operator $I_h$}
\noindent We will denote by $L^2\left(\mathbb{R} \to \mathbb{R}^n \right)$ for $n\in \mathbb{N}$ (or simply $L^2$ where the dimension of the image is understood by context) to be the space of square integrable functions defined on $\mathbb{R}$ with its usual inner product and norm denoted by $\CLIP{\cdot}{\cdot}$ and $\CLN{\cdot}$ respectively (in contrast to $L^2_h$ which the space of square summable functions defined on the mesh $\mathbb{R}_h$ and whose norm is denoted by $\DLN{\cdot}$).  Clearly if we restrict a continuum function $u \in L^2$ to $\mathbb{R}_h$ then we may consider it also as a mesh function in $L^2_h$.  \\
\indent  Lemma 4.1 and corollary 4.2 are based on similar statements found in \cite{Bond, Stum}.  They are the key ingredients allowing us to pass from a discrete function to a continuum function while preserving the necessary estimates for our solution (i.e. boundedness of Schwartz semi-norms).  The proofs of lemma 4.1 and corollary 4.2 can be found in \cite{Bond} where the author uses ideas from \cite{Stum}.\\

\begin{lem}
For any $h>0$ there exists a linear isometry $I_h : L^2_h\left(\mathbb{R}_h \to \mathbb{R}^n \right)  \rightarrow L^2\left(\mathbb{R} \to \mathbb{R}^n \right)$ such that if $u \in L^2_h$ then $U:=I_h u$ has the following properties:
\begin{enumerate}
\item  $U \in C^{\infty}\left(\mathbb{R} \to \mathbb{R}^n \right)$ (hence we can think of $I_h$ as a "smoothing operator").
\item For any point $x_l \in \mathbb{R}_h$ we have that $U (x_l) = u(x_l)$.
\item  For each $j>0$ the following inequalities hold: 
$$\left( \frac{2}{\pi}\right)^j \CLN{\frac{\partial^j}{\partial x^j} U} \leq \DLN{ D^j_+ u} \leq  \CLN{\frac{\partial^j}{\partial x^j} U} $$  
\end{enumerate}
An explicit formula for $U(x)$ is given by

\begin{equation}
U(x)=\sum_{l=-\infty}^{\infty} u(x_l) \left( \frac{ \sin \frac{\pi}{h}(x_l - x)}{\frac{\pi}{h}(x_l-x)} \right) \label{1.1}
\end{equation}
\end{lem}

\begin{cor}
Let $M\geq 2$ be an integer and let $h > 0$ be a real number.  Suppose $u$ is a mesh function on $\mathbb{R}_h$ such that $x_l^N u(x_l) \in L^2_h\left(\mathbb{R}_h\to\mathbb{R}^n \right)$ for each $0 \leq N \leq M$.  Then for each $j \in \mathbb{N}$, and $0 \leq N \leq M-2 $ we have,
$$\left( \frac{2}{\pi}\right)^j \CLN{ \frac{\partial^j}{\partial x^j} \left( x^N U\right) } \leq \DLN{ D^j_+ \left( x^N u\right) } \leq  \CLN{ \frac{\partial^j}{\partial x^j} \left( x^N U \right) } $$  
\end{cor}
\subsection{Schwartz Boundedness of Smoothly Continued Discrete Solutions}
\noindent In this section we will show that a certain smooth continuation of the discrete solution remains bounded in the Schwartz semi-norms $\SupNormtx{\IHN{\cdot}^N \partial_x^n \partial_t^m \cdot} $.  \\
\indent Suppose $u : \mathbb{R}_h \times \mathbb{R}_k\to \mathbb{R}^n $ is a mesh function for some $h,k>0$ which is compactly supported in time for $t_j \in \left[T_0,T_1\right]$ and which satisfies the property that there is some $C>0$ such that $\DLN{u(\cdot,t_j)} \leq C$ for each $t_j \in \mathbb{R}_k$.  We will define a smooth continuation of $u$ by the following:  \\
\indent Since $u$ is compactly supported in time we know that for each $x_l \in \mathbb{R}_h$ we have $u(x_l,\cdot)\in L^2_k$, therefore by lemma 4.1 we may apply the operator $I_k$ to $u(x_l, \cdot)$ in $t$ to obtain by (\ref{1.1}) that for any $t\in \mathbb{R}$,

\begin{equation*}
(I_k u)(x_l,t)  =  \sum_{T_0 \leq t_j \leq T_1 } u(x_l,t_j) \left[ \frac{\sin \frac{\pi}{k}(t-t_j) }{\frac{\pi}{k} (t-t_j)} \right]  
\end{equation*}

\noindent and therefore,

\begin{eqnarray}
\DLN{ I_k u(\cdot,t) } & \leq & \sum_{T_0 \leq t_j \leq T_1 } \left| \frac{\sin \frac{\pi}{k}(t-t_j) }{\frac{\pi}{k} (t-t_j)}  \right| \cdot \DLN{ u(\cdot,t_j)} \nonumber \\
                     & \leq & \sum_{T_0 \leq t_j \leq T_1 } \SupNorm{\frac{\sin y }{y}} \cdot \DLN{ u(\cdot,t_j) } < \infty \label{L10.1}
\end{eqnarray}
 
\noindent Hence by lemma 4.1 we can apply the smoothing operator $I_h$ to $I_k u(\cdot,t)$ in the $x$ variable for each $t\in \mathbb{R}$ to obtain a continuum function $I u:=I_h I_k u$.  By linearity of $I_h$ it follows that $I u$ is given by,
\begin{eqnarray*}
I u(x,t) &=& \sum_{T_0 \leq t_j \leq T_1 } I_h u(x,t_j) \left[ \frac{\sin \frac{\pi}{k}(t-t_j) }{\frac{\pi}{k} (t-t_j)} \right]
\end{eqnarray*}
and since for each $t_j$ we have $\frac{\sin \frac{\pi}{k}(t-t_j)}{\frac{\pi}{k}(t-t_j)}$ is smooth in $t$ and also for each $j$ the function $I_h u(x,t_j)$ is smooth in $x$ we see that $ I u\in C^{\infty}\left(\mathbb{R}\times \mathbb{R} \to \mathbb{R}^n\right)$ and by lemma 4.1 $I u (x_l,t_j) = u(x_l,t_j)$ for any $(x_l,t_j)\in \mathbb{R}_h \times \mathbb{R}_k$.  Moreover, it is clear from the above formula that for any $m \in \mathbb{N}$ we have $\partial_t^m I u = I_h \left( \partial_t^m I_k u \right)$.  Given a discrete solution $u$ from lemma 2.6 we may now construct a smooth continuation $I \hat{u}$ and prove that it remains bounded in the continuum Schwartz seminorms - this is the statement of lemma 4.3 below.  Its proof is the same as the proof of lemma 4.3 in \cite{Gonz} so we will omit it here.  \\
\indent If $u(x,t)=(u_1(x,t),u_2(x,t) ):\mathbb{R} \times \mathbb{R}\to \mathbb{R}^2$ is any function then we will use the notations $
\SupNormx{u(\cdot,t)}:=\SupNorm{u_1(\cdot,t) } + \SupNorm{u_2(\cdot,t)}$ and $\SupNormtx{u}:=\SupNormt{\SupNormx{ u }}$.  

\begin{lem}
Let $u_0 = \left((u_1)_0,(u_2)_0 \right)$ where $(u_l)_0 \in S^{-\infty}(\mathbb{R})$ for $l=1,2$ and let $K,T,\epsilon >0$ all be given as in lemma 2.6 and let $N,n \in \mathbb{N}$, $m \in \left\{0,1,2\right\}$.  Then there exists $0 < k_0 \leq \epsilon$  and $C>0$ both depending on $m,n$ and $N$ such that if $h \in (0,1)$, $k\in (0,k_0)$, and $u(x_n,t_j) = \left( u_1(x_n,t_j) , u_2(x_n,t_j) \right)$ is the solution of the difference scheme (\ref{dgNLS}) defined on $\mathbb{R}_h \times \left( \mathbb{R}_k \cap \left[0,T\right] \right)$ with initial condition $u_0 $ then we have 

\begin{equation*}
\SupNormtx{\IHN{\cdot}^N \partial_x^n \partial_t^m I \hat{u} }\leq C 
\end{equation*}
\end{lem}

\subsection{Proof of Local Existence for the Generalized NLS equation in $S^{-\infty}$}
\noindent By using corollary 4.2 and the Arzela-Ascoli theorem we shall now construct a smooth solution to (\ref{gNLS}) lying in $S^{-\infty}\left(\mathbb{R}\times \left[ 0,T \right]\right)$ that comes from the discrete solution constructed in lemma 2.6.  Theorem 1.1 and its proof are completely analogous to the corresponding results given by Bondareva for the KdV equation (see \cite{Bond} theorem 2) and by the author for the modified KdV equation (see \cite{Gonz} theorem 1.1). \\

\noindent \textbf{Proof of Theorem \ref{thm1.1}} \hspace{1mm} (\textit{existence}) \\  
\indent Since $u_0\in  S^{-\infty}\left(\mathbb{R}\right)$ it follows that there is some $K>0$ such that for any $0<h<1$ we have $\DSN{u_0}\leq K$.  Therefore, by lemma 2.6, there exists $T, L, \epsilon >0$ such that if $h \in (0,1)$ and $k \in (0,\epsilon)$ then there is a solution to the difference scheme (\ref{dgNLS}) with initial condition $u_0$ defined on $\mathbb{R}_h \times \left(\mathbb{R}_k \cap \left[0,T\right]\right)$ and we denote this solution by by $u^{h,k}$.  Let $U^{h,k}=\left( U_1^{h,k},U_2^{h,k} \right) := I \hat{u}^{h,k}$.  From lemma 4.3 we know that for every $N,n \in \mathbb{N}$, $m \in \left\{0,1,2\right\}$ there exists $0 < k_0(m,n,N) \leq \epsilon$ and there exists $C_{m,N,n} > 0$ such that if $h \in (0,1)$, $k\in (0,k_0)$, then we have 

\begin{equation}
\SupNormtx{\IHN{\cdot}^N \partial_x^n \partial_t^m U^{h,k} }\leq C  \label{T.1}
\end{equation}

\vspace{2mm}

\noindent From the family of functions $\left\{U^{h,k}\right\}_{h\in (0,1),k\in (0,k_0(1,1,0))}$ we now wish to extract a convergent subsequence by using the Arzela-Ascoli theorem (this theorem can be found for example in \cite{Rud}).\\
\indent Let $(x_0,t_0),(x_1,t_1)$ be points in $\mathbb{R}\times \mathbb{R}$.  By lemma 4.3 and by the intermediate value theorem we have
\begin{eqnarray}
\left|U^{h,k}(x_0,t_0) - U^{h,k}(x_1,t_1)\right| & \leq & \left|U^{h,k}(x_0,t_0) - U^{h,k}(x_0,t_1)\right| + \left|U^{h,k}(x_0,t_1) - U^{h,k}(x_1,t_1)\right| \nonumber \\
                                                 & \leq & C\left( \left| \partial_t U_1^{h,k}(x_0,\tilde{t}) \right| + \left| \partial_t U_2^{h,k}(x_0,\tilde{\tilde{t}}) \right| \right)\cdot \left|t_0-t_1 \right|  \nonumber \\
                                                 &      &  + C \left( \left| \partial_x U_1^{h,k}(\tilde{x},t_1) \right|+\left| \partial_x U_2^{h,k}(\tilde{\tilde{x}},t_1) \right| \right) \cdot \left|x_0-x_1 \right|  \nonumber \\
                                                 & \leq & \SupNormtx{\partial_t U^{h,k}} \cdot \left| t_0-t_1\right| + \SupNormtx{\partial_x U^{h,k}} \cdot \left|x_0-x_1 \right| \nonumber \\
                                                 & \leq & C \left| (x_0,t_0)-(x_1,t_1) \right|\label{T.2}
\end{eqnarray} 
\noindent which shows that the family of functions $U^{h,k}$ is equicontinuous on $\mathbb{R}\times \mathbb{R}$. \\
\indent From (\ref{T.1}) it follows that the family $U^{h,k}$ is also bounded uniformly for $h\in (0,1)$, \\ $k\in (0,k_0(1,1,0))$.  Hence, by the Arzela-Ascoli theorem we may construct a subsequence $U^{h_i,k_i}$, where of course $h_i,k_i \searrow 0 $ as $i \rightarrow \infty$, converging uniformly on compact sets to a function \\ $U\in C^0 \left( \mathbb{R}\times \mathbb{R} \right) $. \\
\indent The above argument can also be made for the family of functions $\partial_x U^{h_i,k_i} $for $h\in (0,1)$, \\ $k\in (0,k_0(1,2,0))$.  Namely, estimate (\ref{T.1}) implies that the family is bounded uniformly and also that we may use estimate (\ref{T.2}) with $U^{h,k} $ replaced by $\partial_x U^{h_i,k_i}$ to see that it is also an equicontinuous family.  Thus we conclude that there is some $V\in C^0\left( \mathbb{R}\times \mathbb{R} \right)$ and a subsequence  $\partial_x U^{h_l,k_l}$ converging uniformly on compact sets to $V$.  Since we have uniform convergence on compact sets for $U^{h_l, k_l}$ and $\partial_x U^{h_l,k_l}$ it follows that $U$ is differentiable in $x$ and $\partial_x U = V$ on $\mathbb{R}\times\mathbb{R}$.  \\
\indent By repeating the same argument we conclude by induction that for each $p\in\mathbb{N}$ the function $\partial_x^{p-1} U\in C^0(\mathbb{R}\times\mathbb{R})$ is differentiable in $x$ because the sequence $\partial_x^p U^{h_l,k_l}$ for $h\in (0,1)$ and \\ $k\in (0,k_0(1,p+1,0))$ is bounded uniformly by (\ref{T.1}) and is equicontinuous by (\ref{T.2}) and hence it has a subsequence uniformly convergent on compact subsets of $\mathbb{R}\times \mathbb{R}$ to $\partial_x \partial_x^{p-1} U$.  In this way we will obtain a countable array of subsequences, one for each $p\in \mathbb{N}$ and from this array we extract a diagonal subsequence.  From this diagonal subsequence it will follow that for each $p\in \mathbb{N}$ we have $\partial_x^p U^{h_l,k_l} \rightarrow \partial_p U$ uniformly on compact sets.\\
\indent Consider the family of functions $\partial_t U^{h_l,k_l}$ for $h\in (0,1)$ and $k\in (0,k_0(2,1,0))$.  Estimate (\ref{T.1}) for $m=1$ implies that the family is bounded uniformly and also that we may use estimate (\ref{T.2}) with $U^{h,k}$ replaced by $\partial_t U^{h_i,k_i}$ to see that it is also an equicontinuous family.  Hence we may, as before for $x$, conclude that $U$ is differentiable in $t$ and construct a subsequence of $\partial_t U^{h_l,k_l}$ uniformly convergent on compact sets to $\partial_t U$.  From this subsequence of $(h_l,k_l)$ we consider the family $\partial_x \partial_t U^{h_l,k_l}$ for $h\in (0,1)$ and $k\in (0,k_0(2,2,0))$.  Again estimate (\ref{T.1}) implies that the family is bounded uniformly and also we may use estimate (\ref{T.2}) with $U^{h,k}$ replaced by $\partial_x\partial_t U^{h_l,k_l}$ to see that it is also an equicontinuous family.   Thus we may again extract a subsequence $U^{h_l,k_l}$ to see that $\partial_t U$ is differentiable in $x$ and $\partial_x \partial_t U^{h_l,k_l}\rightarrow \partial_x \partial_t U$ uniformly on compact sets.   Continuing inductively we consider the sequence of functions $\partial_x^p \partial_t U^{h_l,k_l}$ for $h\in (0,1)$, $k\in (0,k_0(2,p+1,0))$.  It is equicontinuous by (\ref{T.2}) and from (\ref{T.1}) it is uniformly bounded, thus we conclude that $\partial_x^{p-1}\partial_t U$ is differentiable in $x$ and we may extract a subsequence so that $\partial_x^p \partial_t U^{h_l,k_l}\rightarrow \partial_x^p \partial_t U$.  \\
\indent Continuing in this way we will again obtain an array of subsequences of $U^{h_l,k_l}$, one for each $p\in \mathbb{N}$.  By taking a diagonal subsequence we obtain a subsequence such that for each $p\in \mathbb{N}$ and for $q=0,1$ we have $\partial_x^p\partial_t^q U^{h_l,k_l} \rightarrow \partial_x^p\partial_t^q U$ uniformly on compact subsets of $\mathbb{R}\times \mathbb{R}$.  In addition, it follows that for any $N\in \mathbb{N}$ the sequence $\IHN{x}^N \partial_x^p \partial_t^q U^{h_l,k_l} \rightarrow \IHN{x}^N \partial_x^p\partial_t^q U$ uniformly on compact subsets of $\mathbb{R}\times \mathbb{R}$ because for any compact set $X \subset \mathbb{R}\times\mathbb{R}$ we have the inequality
\begin{eqnarray*}
\left|\IHN{x}^N \partial_x^p \partial_t^q U^{h,k}(x,t) - \IHN{x}^N \partial_x^p \partial_t^q U(x,t)\right| & \leq & \max_{x\in X}\IHN{x}^N \cdot \left| \partial_x^p \partial_t^q U^{h,k}(x,t) - \partial_x^p\partial_t^q U(x,t)\right|
\end{eqnarray*} 

\noindent By construction we can see that $U$ satisfies the following conditions:\\
\begin{enumerate}
\item $\partial_x^p U$ exists for each $p\in \mathbb{N}$ and is continuous (i.e. it lies in $C^0(\mathbb{R}\times\mathbb{R})$).\\
\item $\partial_x^q \partial_t \partial_x^p U$ exists for each $p,q\in \mathbb{N}$ and is continuous (i.e. it lies in $C^0(\mathbb{R}\times\mathbb{R})$).\\
\item If $p+q = p'+q'$ then $\partial_x^q \partial_t \partial_x^p U = \partial_x^{q'} \partial_t \partial_x^{p'} U$. 
\end{enumerate}

\indent We now claim that the complex-valued function $U_1+i U_2$ is a solution to (\ref{gNLS}) (and we shall also identify this solution by $U$).  To prove it we fix a point $(x,t) \in \mathbb{R} \times \left[0,T\right]$ and show that the system (\ref{rgNLS}), (\ref{igNLS}) is satisfied by its components $U_1$ and $U_2$ at $(x,t)$.  From our final subsequence of pairs $(h_l,k_l)$ above we first construct points $(x_l,t_l) \in \mathbb{R}_{h_l}\times\left( \mathbb{R}_{k_l}\cap \left[0,T\right] \right)$ to be the nearest points in the grid to $(x,t)$ (note: in this context $x_l \neq l \cdot h_l$ and $t_l \neq l \cdot k_l$).  It then follows that $(x_l,t_l) \rightarrow (x,t)$ as $l\rightarrow \infty$.  By construction we have for each $l\in \mathbb{N}$ that $\left(U_1^{h_l,k_l},U_2^{h_l,k_l} \right)$ satisfies the pair of difference equations (\ref{drgNLS}), (\ref{digNLS}) at point $(x_l,t_l)$.  Replace the discrete derivatives in $t$ and $x$ of equations (\ref{drgNLS}) and (\ref{digNLS}) by ordinary derivatives at intermediate points $(\tilde{x}_l,\tilde{t}_l)$ (possibly different intermediate points for each term containing derivatives).  We will then obtain a sum of products of terms of the form  $\partial_t U_r^{h_l,k_l}(x_l,\tilde{t}_l)$, $U_r^{h_l,k_l}(x_l,t_l)$, $U_r^{h_l,k_l}(x_l,t_l +k_l)$, $\partial_x^2 U_r^{h_l,k_l}(\tilde{x}_l,t_l)$, $f_r(x_l,t_l)$, and $g_r(x_l,t_l)$ for $r=1,2$.  By continuity of $f$ and $g$ we see that $f(x_l,t_l) \rightarrow f(x,t)$ and $g(x_l,t_l)\rightarrow g(x,t)$ as $l\rightarrow \infty$.  Moreover, since $\left(U_1^{h_l,k_l}(x,t),U_2^{h_l,k_l}(x,t) \right) \rightarrow \left( U_1(x,t),U_2(x,t) \right)$ as $l \rightarrow \infty$ and
\begin{align*}
& \left|\left(U_1^{h_l,k_l}(x_l,t_l),U_2^{h_l,k_l}(x_l,t_l)\right)-\left( U_1(x,t),U_2(x,t) \right) \right| \\
& \leq  \left|\left(U_1^{h_l,k_l}(x_l,t_l),U_2^{h_l,k_l}(x_l,t_l)\right) - \left(U_1^{h_l,k_l}(x,t),U_2^{h_l,k_l}(x,t)\right)\right| \\
&   + \left|\left(U_1^{h_l,k_l}(x,t),U_2^{h_l,k_l}(x,t)\right) - \left( U_1(x,t),U_2(x,t) \right) \right|
\end{align*}

\noindent it follows by equicontinuity of the family $\left(U_1^{h_l,k_l},U_2^{h_l,k_l}\right)$ for $h,k$ sufficiently small that \\ $\left(U_1^{h_l,k_l},U_2^{h_l,k_l}\right) \rightarrow \left( U_1(x,t),U_2(x,t) \right)$.  We may use the same convergence argument for the other terms in the equation to show that as $l\rightarrow \infty$ the equations (\ref{drgNLS}), (\ref{digNLS}) becomes equations (\ref{rgNLS}),(\ref{igNLS}) at the point $(x,t)$.  Therefore $\left(U_1, U_2\right)$ satisfies the equations (\ref{rgNLS}),(\ref{igNLS}) and hence $U = U_1+iU_2$ satisfies (\ref{gNLS}).\\ 
\indent Since $U$ satisfies the equation (\ref{gNLS}) it follows that $\partial_t U$ is also differentiable in time and its higher time derivatives can be written in terms of the lower $x$ derivatives.  The derivatives also clearly commute as was mentioned above in condition 3, therefore it follows that $U\in C^{\infty}\left(\mathbb{R}\times \left[0,T\right] \to \mathbb{C} \right)$.\\ 
\indent Furthermore, we can show that the limit function $U$ is in $S^{-\infty}\left(\mathbb{R}\times\left[0,T\right] \to \mathbb{C} \right)$.  By taking the limit of $\IHN{x}^N \partial_x^n \partial_t^m U^{h_l,k_l}(x,t)$ as $l \rightarrow \infty$ we can see that (\ref{T.1}) also holds for the function $U$ with $n,N \in \mathbb{N}$ and $m=0,1$.  By repeatedly using the equation (\ref{gNLS}) we may write $\IHN{\cdot}^N \partial_x^n \partial_t^m U$ as a sum of products of terms of the form $\IHN{\cdot}^N \partial_x^n U$ each of which can be bounded by some constant depending on $N,n\in\mathbb{N}$ by using the limiting case of (\ref{T.1}) and this implies that for any $m,n,N\in \mathbb{N}$ we have, \\ 
\begin{eqnarray*}
\SupNormtx{\IHN{\cdot}^N \partial_x^n \partial_t^m U}  &\leq&  C_{m,n,N} 
\end{eqnarray*}
\noindent which shows that $U\in S^{-\infty}\left(\mathbb{R}\times\left[0,T\right]\to \mathbb{C}\right)$. \hspace{4cm} $\square$ 

\subsection{Proof of Local Existence for the NLS equation in $S^{\beta}$ when $\beta \leq 0$}
\noindent Now we shall now construct smooth solutions to (\ref{NLS}) lying in $S^{\beta}\left(\mathbb{R}\times \left[ 0,T \right] \to \mathbb{C}\right)$ for $\beta \leq 0$ that come from adding an above solution of (\ref{gNLS}) to the asymptotic solution constructed in lemma \ref{lemA.2}. \\

\noindent \textbf{Proof of Theorem \ref{thm1.2}} \hspace{1mm} (\textit{existence}) \\
\indent  By lemma \ref{lemA.2} there exists an asymptotic solution $f(x,t) \in S^{\beta}(\mathbb{R}\times \mathbb{R}\to \mathbb{C})$ of the initial value problem (\ref{NLS}) whose expansion coefficients satisfy the desired property.  Let $u_0(x) = w_0(x)-f(x,0)$ and let $g: = if_t+ f_{xx}+\mu \overline{f}f^2 $.   By construction $u_0 \in S^{-\infty}(\mathbb{R}\to \mathbb{C})$.  Moreover $f$ and $g$ satisfy the hypotheses of theorem 1.1. Therefore there exists a $T>0$ and a solution $u(x,t) \in S^{-\infty}(\mathbb{R}\times \left[0,T\right])$ to equation (\ref{gNLS}).  Let $w(x,t):=u(x,t)+f(x,t)$.  Since $u$ satisfies (\ref{gNLS}) it follows that $w$ satisfies (\ref{NLS}). 
 Moreover, since $u \sim 0$ it follows that $w$ and $f$ have the same asymptotic expansions and in particular the coefficients in the asymptotic expansions of $w$ satisfy the second statement of the theorem.  Finally, since $f(x,t) \in S^{\beta}(\mathbb{R}\times \left[0,T\right] \to \mathbb{C})$ it follows that $w(x,t) \in S^{\beta}(\mathbb{R}\times \left[0,T\right] \to \mathbb{C})$.  \hspace{2cm}  $\square$ \\
\section{Uniqueness of Solutions}
\subsection{Uniqueness in $S^{-\infty}$ for the Generalized NLS Equation}
\noindent In this section we shall prove uniqueness of solutions in $S^{-\infty}(\mathbb{R}\times\left[0,T\right])$ for (\ref{gNLS}) by using Gronwall's Inequality.  We shall state this inequality below and we refer the reader to \cite{Gonz} for a short proof. 
\begin{lem} \label{Gron}
\noindent Let $T>0$ and $c_1,c_2 \in \mathbb{R}$ be given and $c_1 \neq 0$.  Suppose $\eta: \left[0,T\right] \rightarrow \mathbb{R}$ is a nonnegative, differentiable function and that for each $t\in \left[0,T\right]$ we have\\
\begin{equation*}
\frac{d\eta}{dt}(t)  \leq  c_1 \eta(t) + c_2
\end{equation*}
\noindent Then for each $t\in \left[0,T\right]$ we have\\
\begin{equation*}
\eta(t)  \leq  e^{c_1 t}\left(\eta(0)-\frac{c_1}{c_2}\right) - \frac{c_1}{c_2}
\end{equation*}
\end{lem}

\noindent \textbf{Proof of Theorem \ref{thm1.1}} \hspace{1mm} (\textit{uniqueness}) \\
\indent  Suppose $u(x,t)=u_1(x,t)+iu_2(x,t),v(x,t)=v_1(x,t)+iv_2(x,t) \in S^{-\infty}(\mathbb{R}\times\left[0,T\right] \to \mathbb{R})$ are two solutions of (\ref{gNLS}) with initial data $u_0 = \left( (u_1)_0,(u_2)_0 \right) \in S^{-\infty}(\mathbb{R}\to \mathbb{C})$.  Then by equations (rgNLS) and (igNLS) we have
\begin{eqnarray*}
(u_1)_t + (u_2)_{xx} + \mu \left( u_1^2 u_2 + u_2^3 + 2u_1 u_2 f_1 + u_1^2 f_2 + 3 u_2^2 f_2 + 2f_1 f_2 u_1 + u_2 f_1^2 + 3u_2 f_2^2\right) + g_2 & = &  0 \\
(v_1)_t + (v_2)_{xx} + \mu \left( v_1^2 v_2 + v_2^3 + 2v_1 v_2 f_1 + v_1^2 f_2 + 3 v_2^2 f_2 + 2f_1 f_2 v_1 + v_2 f_1^2 + 3v_2 f_2^2\right) + g_2 & = &  0
\end{eqnarray*}
\noindent and
\begin{eqnarray*}
(u_2)_t - (u_1)_{xx} - \mu \left( u_1^3 + u_1 u_2^2 + u_1 f_2^2 + 3u_1^2 f_1 + u_2^2 f_1 + 2u_1 u_2 f_2 + 3 u_1 f_1^2 + 2 f_1 f_2 u_2 \right) - g_1  & = &  0 \\
(v_2)_t - (v_1)_{xx} - \mu \left( v_1^3 + v_1 v_2^2 + v_1 f_2^2 + 3v_1^2 f_1 + v_2^2 f_1 + 2v_1 v_2 f_2 + 3 v_1 f_1^2 + 2 f_1 f_2 v_2 \right) - g_1  & = &  0
\end{eqnarray*}

\noindent Let $q(x,t) = (q_1(x,t),q_2(x,t)) := u(x,t)-v(x,t)$.  By subtracting the above two equations we see that the components of $q$ satisfy the equations
\begin{eqnarray*}
(q_1)_t + (q_2)_{xx} + \mu \left[ u_1 u_2 q_1 + u_1 v_1 q_2 + v_1 v_2 q_1 + q_2 \left( u_2^2 + u_2 v_2 + v_2^2\right) \right. \hspace{3.5cm}\\ \left. + 2f_1 u_2 q_1 + 2f_1 v_1 q_2 + f_2 q_1 \left( u_1+v_1 \right) + 3 f_2 q_2 \left( u_2+v_2\right) + 2f_1 f_2 q_1 + f_1^2 q_2 + 3f_2^2 q_2 \right]  =  0\\
(q_2)_t - (q_1)_{xx} - \mu \left[ u_1 u_2 q_2 + u_1 v_2 q_2 + v_2^2 q_1 + q_1 \left( u_1^2+u_1 v_1 + v_2^2\right) \right. \hspace{3.5cm} \\ \left.
+ f_2^2q_1 + 3f_1 q_1 \left( u_1+v_1\right)+ f_2 q_2 \left( u_2+v_2 \right) + 3f_1^2 q_1 + 2f_1 f_2 q_2 + 2f_2 u_2 q_1 + 2 f_1 v_1 q_2 \right]  =  0
\end{eqnarray*}
\noindent Thus if we multiply the first equation by $q_1$ and the second equation by $q_2$ then add them together and integrate by parts in $x$ over $(-\infty,\infty)$ we get an estimate of the form
\begin{equation*}
\frac{d}{dt} \int_{\mathbb{R}} \left(q_1^2 + q_2^2 \right) dx \leq  C_{u,v,f,T} \int_{\mathbb{R}} \left( q_1^2 + q_2^2\right) dx
\end{equation*}
\noindent and moreover, $q(\cdot,0) = 0$, therefore by lemma 5.1 it follows that $\CLN{q(\cdot,t)} = 0$ for all $t\in \left[0,T\right]$ and since $q$ is smooth this implies that $q(x,t)=0$ for all $(x,t) \in \mathbb{R}\times\left[0,T\right]$. \hspace{2cm} $\square$ 

\subsection{Uniqueness in $S^{\beta}$ for the NLS Equation When $\beta \leq 0$}
\noindent In this section we shall prove uniqueness of solutions in $S^{\beta}(\mathbb{R}\times\left[0,T\right])$ for (\ref{NLS}) when $\beta \leq 0$.  First we will need the following lemma. 
\begin{lem} \label{lem5.1}
Let $I\subset \mathbb{R}$ be an interval and $\beta \leq 0$.  Suppose $w(x,t) \in S^{\beta}(\mathbb{R}\times I \to \mathbb{C})$ is a solution to (\ref{NLS}) with initial data $w_0\in S^{\beta}(\mathbb{R} \to \mathbb{C})$ and that $ w(x,t) \sim \sum_{k=0}^{\infty} \left(a_k^{\pm}(t) + i b_k^{\pm}(t) \right)x^{\beta_k}$ as $x\to \pm \infty$.  Then $ \sum_{k=0}^{\infty} \left(a_k^{\pm}(t) + i b_k^{\pm}(t) \right)x^{\beta_k}$ is a formal solution to (\ref{NLS}). 
\end{lem}

\noindent \textbf{Proof of Lemma \ref{lem5.1}} \hspace{1mm} By symmetry it suffices to show that the positive $x$ asymptotic expansion satisfies equation (\ref{NLS}).   Let $A_0=\{\beta_j \}_{j=0}^{\infty}$ and let $J \subset I$ be a compact interval.  We enlarge $A_0$ to the set $\Gamma$ defined in appendix A having the properties mentioned in lemma \ref{lemA.1}.  Let us re-write the asymptotic expansion as $\sum_{k=0}^{\infty}\left(  a_k^+(t) + i b_k^{+}(t)\right) x^{\gamma_k}$ where $a_k^+(t)+ib_k^{+}(t)=0$ if $\gamma_k \notin A_0$.  By definition of being asymptotic it follows that for every $N \in \mathbb{N}$ we may write $$w(x,t)=\sum_{k=0}^N \left(  a_k^+(t) + i b_k^{+}(t)\right) x^{\gamma_k} + R_N(x,t)$$ for $x>1$ and $t\in J$ where $\partial_t^i \partial_x^j R_N(x,t) = O\left(|x|^{\gamma_{N+1}-j} \right)$ for every $i,j \in \mathbb{N}$.  Let \\ $f_N(x,t) = \sum_{k=0}^N \left(  a_k^+(t) + i b_k^{+}(t)\right)x^{\gamma_k}$.  Then by plugging $w$ into (\ref{NLS}) and separating real and imaginary parts, it follows as we see in (\ref{a.1}) and (\ref{a.2}) that for some $M \leq N$ we have
\begin{eqnarray}
\sum_{j=0}^M \left[ \dot{a}_j^+ + \mu \sum_{{l,m,n}\atop {\gamma_l+\gamma_m+\gamma_n=\gamma_j}} - a_l^+ a_m^+  b_n^+ + a_l^+ a_n^+  b_m^+ +a_n^+  a_m^+  b_l^+ +b_l^+  b_m^+  b_n^+ + \sum_{{p}\atop{\gamma_p-2=\gamma_j}} b_p^+  \gamma_p  (\gamma_p - 1) \right]x^{\gamma_j}\nonumber \\ + O\left(|x|^{2\gamma_0+\gamma_{N+1}}\right) =0 \hspace{2cm} \label{L12.11} \hspace{1cm}
\end{eqnarray}
and 
\begin{eqnarray}
\sum_{j=0}^M \left[ \dot{b}_j^+ -  \mu \sum_{{l,m,n}\atop {\gamma_l+\gamma_m+\gamma_n=\gamma_j}} a_l^+  a_m^+  a_n^+ + a_l^+  b_n^+  b_m^+ +b_n^+  a_m^+  b_l^+ - b_l^+  b_m^+  a_n^+  - \sum_{{p}\atop{\gamma_p-2=\gamma_j}} a_p^+  \gamma_p (\gamma_p - 1) \right] x^{\gamma_j} \nonumber \\ + O\left(|x|^{2\gamma_0+\gamma_{N+1}}\right) =0 \hspace{2cm} \label{L12.12} \hspace{1cm}
\end{eqnarray}

\noindent We may assume that $N$ is sufficiently large so that $M \geq 1$ and $2\gamma_0+\gamma_{N+1} < \gamma_1$.  Since the above equation must hold for all $x>1$ we may divide by $x^{\gamma_0}$ to obtain from (\ref{L12.11}) and (\ref{L12.12}) that 

\begin{align*}
\dot{a}_0^+ + & \mu \sum_{{l,m,n}\atop {\gamma_l+\gamma_m+\gamma_n=\gamma_0}} - a_l^+ a_m^+  b_n^+ + a_l^+ a_n^+  b_m^+ +a_n^+  a_m^+ b_l^+ +b_l^+  b_m^+  b_n^+ \\
              & + \sum_{{p}\atop{\gamma_p-2=\gamma_0}} b_p^+ \gamma_p (\gamma_p - 1) + O(|x|^{\gamma_1-\gamma_0}) =  0
\end{align*}
and
\begin{align*}
\dot{b}_0^+ -  & \mu \sum_{{l,m,n}\atop {\gamma_l+\gamma_m+\gamma_n=\gamma_0}} a_l^+  a_m^+ a_n^+ + a_l^+  b_n^+  b_m^+ +b_n^+  a_m^+  b_l^+ - b_l^+ b_m^+  a_n^+  \\
               & - \sum_{{p}\atop{\gamma_p-2=\gamma_0}} a_p^+ \gamma_p (\gamma_p - 1)  + O(|x|^{\gamma_1-\gamma_0}) =  0
\end{align*}

\noindent and hence 
\begin{eqnarray*}
\dot{a}_0^+ + \mu \sum_{{l,m,n}\atop {\gamma_l+\gamma_m+\gamma_n=\gamma_0}} - a_l^+ a_m^+  b_n^+ + a_l^+ a_n^+  b_m^+ +a_n^+  a_m^+  b_l^+ +b_l^+  b_m^+  b_n^+ + \sum_{{p}\atop{\gamma_p-2=\gamma_0}} b_p^+  \gamma_p (\gamma_p - 1)  =  0
\end{eqnarray*}
and
\begin{eqnarray*}
\dot{b}_0^+ -  \mu \sum_{{l,m,n}\atop {\gamma_l+\gamma_m+\gamma_n=\gamma_0}} a_l^+  a_m^+  a_n^+ + a_l^+  b_n^+  b_m^+ +b_n^+  a_m^+  b_l^+ - b_l^+  b_m^+  a_n^+  - \sum_{{p}\atop{\gamma_p-2=\gamma_0}} a_p^+  \gamma_p  (\gamma_p - 1) =  0
\end{eqnarray*}
\indent Continuing in the same way we may assume that $N$ is sufficiently large so that $2\gamma_0+\gamma_{N+1}-1 < \gamma_2$.  Dividing (\ref{L12.11}) and (\ref{L12.12}) by $x^{\gamma_1}$ we obtain that 
\begin{align*}
\dot{a}_1^+ + & \mu \sum_{{l,m,n}\atop {\gamma_l+\gamma_m+\gamma_n=\gamma_1}} - a_l^+ a_m^+  b_n^+ + a_l^+ a_n^+  b_m^+ +a_n^+  a_m^+ b_l^+ +b_l^+  b_m^+  b_n^+ \\
              & + \sum_{{p}\atop{\gamma_p-2=\gamma_1}} b_p^+ \gamma_p (\gamma_p - 1) + O(|x|^{\gamma_2-\gamma_1}) =  0
\end{align*}
and
\begin{align*}
\dot{b}_1^+ -  & \mu \sum_{{l,m,n}\atop {\gamma_l+\gamma_m+\gamma_n=\gamma_1}} a_l^+  a_m^+ a_n^+ + a_l^+  b_n^+  b_m^+ +b_n^+  a_m^+  b_l^+ - b_l^+ b_m^+  a_n^+  \\
               & - \sum_{{p}\atop{\gamma_p-2=\gamma_1}} a_p^+ \gamma_p (\gamma_p - 1)  + O(|x|^{\gamma_2-\gamma_1}) =  0
\end{align*}

\noindent and hence 
\begin{eqnarray*}
\dot{a}_1^+ + \mu \sum_{{l,m,n}\atop {\gamma_l+\gamma_m+\gamma_n=\gamma_1}} - a_l^+ a_m^+  b_n^+ + a_l^+ a_n^+  b_m^+ +a_n^+  a_m^+  b_l^+ +b_l^+  b_m^+  b_n^+ + \sum_{{p}\atop{\gamma_p-2=\gamma_1}} b_p^+  \gamma_p (\gamma_p - 1)  =  0
\end{eqnarray*}
and
\begin{eqnarray*}
\dot{b}_1^+ -  \mu \sum_{{l,m,n}\atop {\gamma_l+\gamma_m+\gamma_n=\gamma_1}} a_l^+  a_m^+  a_n^+ + a_l^+  b_n^+  b_m^+ +b_n^+  a_m^+  b_l^+ - b_l^+  b_m^+  a_n^+  - \sum_{{p}\atop{\gamma_p-2=\gamma_1}} a_p^+  \gamma_p  (\gamma_p - 1) =  0
\end{eqnarray*}
\indent This process may be repeated inductively to obtain from (\ref{L12.11}) and (\ref{L12.12}) that for any $j\in \mathbb{N}$ we have $$ \dot{a}_j^+ + \mu \sum_{{l,m,n}\atop {\gamma_l+\gamma_m+\gamma_n=\gamma_j}} - a_l^+ a_m^+  b_n^+ + a_l^+ a_n^+  b_m^+ +a_n^+  a_m^+  b_l^+ +b_l^+  b_m^+  b_n^+ + \sum_{{p}\atop{\gamma_p-2=\gamma_j}} b_p^+  \gamma_p  (\gamma_p - 1) =0$$ and $$\dot{b}_j^+ -  \mu \sum_{{l,m,n}\atop {\gamma_l+\gamma_m+\gamma_n=\gamma_j}} a_l^+  a_m^+  a_n^+ + a_l^+  b_n^+  b_m^+ +b_n^+  a_m^+  b_l^+ - b_l^+  b_m^+  a_n^+  - \sum_{{p}\atop{\gamma_p-2=\gamma_j}} a_p^+  \gamma_p (\gamma_p - 1) =0$$  and hence $\sum_{k=0}^{\infty} \left( a_k^+(t) +i b_k^+(t) \right)x^{\gamma_k}$ is a formal solution to (\ref{NLS}). \hspace{2cm} $\square$ \\

\noindent \textbf{Proof of Theorem \ref{thm1.2}} \hspace{1mm} (\textit{uniqueness}) \\
\indent Suppose $w(x,t),r(x,t)\in S^{\beta}(\mathbb{R}\times \left[0,T\right] \to \mathbb{C})$ are two solutions of (\ref{NLS}) with initial data \\ $w_0(x)\in S^{\beta}(\mathbb{R} \to \mathbb{C})$ and that 

$$w_0(x) \sim \sum_{k=0}^{\infty} \left( p_k^{\pm} + i q_k^{\pm} \right) x^{\beta_k} \hspace{.5cm} w(x,t) \sim \sum_{k=0}^{\infty} \left( a_k^{\pm}(t) +i b_k^{\pm}(t) \right) x^{\alpha_k} \hspace{.5cm} r(x,t) \sim \sum_{k=0}^{\infty} \left(c_k^{\pm}(t) + i d_k^{\pm}(t)\right) x^{\delta_k}$$ as $x\to \pm \infty$.  Let $B_0=\left\{\beta_k\right\}_{k=0}^{\infty}$, $A_0=\left\{\alpha_k\right\}_{k=0}^{\infty}$, and $D_0=\left\{\delta_k\right\}_{k=0}^{\infty}$.  By lemma \ref{lem5.1} \\ $\sum_{k=0}^{\infty} \left( a_k^{\pm}(t) +i b_k^{\pm}(t) \right) x^{\alpha_k}$ and $\sum_{k=0}^{\infty} \left(c_k^{\pm}(t) + i d_k^{\pm}(t)\right) x^{\delta_k}$ are formal solutions to (\ref{NLS}) with initial data $\sum_{k=0}^{\infty} \left( p_k^{\pm} + i q_k^{\pm} \right) x^{\beta_k}$ and hence we may assume that $B_0 \subset A_0$ and $B_0 \subset D_0$.  Let $\Lambda = A_0 \cup D_0$ and $\Gamma=\left\{\gamma_k\right\}_{k=0}^{\infty}$ be the set constructed in appendix A from $\Lambda$ having the properties stated in lemma \ref{lemA.1}.  Then after reindexing we may rewrite the asymptotic expansions for $w_0$, $w(x,t)$, and $r(x,t)$ as  $$w_0(x) \sim \sum_{k=0}^{\infty} \left( p_k^{\pm} + i q_k^{\pm} \right) x^{\gamma_k} \hspace{.5cm} w(x,t) \sim \sum_{k=0}^{\infty} \left( a_k^{\pm}(t) +i b_k^{\pm}(t) \right) x^{\gamma_k} \hspace{.5cm} r(x,t) \sim \sum_{k=0}^{\infty}\left(c_k^{\pm}(t) + i d_k^{\pm}(t)\right) x^{\gamma_k}$$ where $p_k^{\pm} + i q_k^{\pm} =0$ if $\gamma_k \notin B_0$, $a_k^{\pm}(t) +i b_k^{\pm}(t)=0$ if $\gamma_k\notin A_0$, and $c_k^{\pm}(t) + i d_k^{\pm}(t) = 0$ if $\gamma_k \notin D_0$.  Since $\sum_{k=0}^{\infty} \left( a_k^{\pm}(t) +i b_k^{\pm}(t) \right) x^{\gamma_k}$ and $\sum_{k=0}^{\infty} \left(c_k^{\pm}(t) + i d_k^{\pm}(t)\right) x^{\gamma_k}$ are formal solutions with initial data $\sum_{k=0}^{\infty} \left( p_k^{\pm} + i q_k^{\pm} \right) x^{\gamma_k}$ it follows that the coefficient pairs $\left( a_k^{\pm}, b_k^{\pm} \right)$ and $ \left( c_k^{\pm}, d_k^{\pm} \right)$ both satisfy equations (\ref{a.1}) and (\ref{a.2}) with the same initial data.  Hence for all $k\in \mathbb{N}$  and all $t\in \left[0,T\right]$ we have $a_k^{\pm}(t) +i b_k^{\pm}(t) = c_k^{\pm}(t) + i d_k^{\pm}(t)$ so that $w(x,t)-r(x,t) \in S^{-\infty}(\mathbb{R} \times \left[0,T\right] \to \mathbb{C})$.  \\
\indent Let $u(x,t) = w(x,t)-r(x,t)$.  Then $u(x,t)$ satisfies (\ref{gNLS}) with initial condition $u_0(x)=0$ and where $f(x,t) = r(x,t)$ and $g(x,t) = 0$.  By uniqueness of solutions to (\ref{gNLS}) in $S^{-\infty}(\mathbb{R}\times \left[0,T\right])$, which was proven in theorem 1.1, it follows that $u(x,t) = 0$ for all $(x,t) \in \mathbb{R}\times \left[0,T\right]$. \hspace{2cm} $\square$ \\

\noindent \textbf{ \large{Acknowledgements}} \\

\indent  The author wishes to thank P. Topalov for suggesting this problem and J. Holmer for kindly providing some related references.  
\appendix
\section{Appendix:  Existence of an Asymptotic Solution}
\noindent In this section we will prove existence of an asymptotic solution to (\ref{NLS}).  First we will need the following lemma.  \\
\indent  Let $A_0=\{\beta_j \}_{j=0}^{\infty}$.  Where $0 \geq \beta_0> \beta_1 > \cdots$, and $\lim_{j\to \infty}\beta_j =-\infty$.  We enlarge $A_0$ to the set $\Gamma$ given by,
\begin{equation*}
\Gamma := \left\{ \sum_{p=1}^k\beta_{i_p} - 2l : k\geq 1, l\geq 0, k,l \in \mathbb{Z}, \beta_{i_p}\in A_0\right\}
\end{equation*}

\vspace{2mm}

\begin{lem} \label{lemA.1}
The set $\Gamma$ has the following properties:
\begin{enumerate}
\item $A_0 \subset \Gamma$
\item $\Gamma$ is countable.
\item $\Gamma$ is bounded above by $\beta_0$.
\item If $\gamma_l$, $\gamma_m$, $\gamma_n$ are all in $\Gamma$ then $\gamma_l+\gamma_m+\gamma_n $ is in $\Gamma$.
\item If $\gamma_p$ is in $\Gamma$ then $\gamma_p - 2$ is also in $\Gamma$.
\item $\Gamma$ is lower finite,  i.e. $\Gamma \cap \left[\left.-M,\infty\right.\right)$ is finite for every $M>0$.
\end{enumerate}
\end{lem}

\vspace{2mm}

\noindent \textbf{Proof of Lemma \ref{lemA.1}} \hspace{1mm} Statements 1 to 5 follow easily from the definition of $\Gamma$ so we shall only prove lower finiteness here.\\
\indent  Let $\Sigma A :=\left\{ \sum_{p=1}^k\beta_{i_p} : k\geq 1, k \in \mathbb{Z}, \beta_{i_p}\in A_0\right\}$.  Since $A_0$ is lower finite it follows that $\Sigma A$ is also lower finite.  Therefore $\Sigma A-\Big(2\mathbb{N} \Big) = \Gamma $ is lower finite. $\square$ \\
\begin{lem} \label{lemA.2}
For any $\beta \leq 0 $ and for any initial condition $w_0 \in S^{\beta}(\mathbb{R})$ there exists an asymptotic solution $f(x,t) \in S^{\beta}(\mathbb{R}\times \mathbb{R})$ of the initial value problem (\ref{NLS}).  Moreover, if $w_0 \sim \sum_{k=0}^{\infty} \left( a_k^{\pm} +i b_k^{\pm}\right)x^{\beta_k}$ and $j$ is the smallest index such that $a_j^{+} + i b_j^{+} \ne 0$ (resp. $a_j^{-}+i b_j^{-} \ne 0$) then the coefficient $a_j^{+}(t) + i b_j^{+}(t)$ (resp. $a_j^{-}(t)+i b_j^{-}(t)$) in the asymptotic expansion of the solution is a nonvanishing continuous function of $t$ and all preceeding coefficients are identically zero.
\end{lem}
\noindent \textbf{Proof of Lemma \ref{lemA.2}} \hspace{1mm} First we will show how to construct a formal solution $\sum_{k=0}^{\infty} \left( a_k^{\pm} +i b_k^{\pm}\right) x^{\beta_k}$.
By symmetry it suffices to construct only the positive $x$ formal solution.  For simplicity we shall omit the superscript $+$ sign in the coefficients $a_j(t)$. \\
\indent First we enlarge the exponent set $A_0=\{\beta_j \}_{j=0}^{\infty}$ to the set $\Gamma$ as defined above.  From lemma \ref{lemA.1} it follows that we may write the set $\Gamma$ as a decreasing sequence $\Gamma=\left\{\gamma_j\right\}_{j=0}^{\infty}$ where $0 \geq \gamma_0$, $\gamma_j >\gamma_{j+1}$, and $\gamma_j \rightarrow -\infty$ as $j\rightarrow \infty$, and we may rewrite the positive $x$ asymptotic expansion of $w_0$ as $\sum_{j=0}^{\infty}(a_j +i b_j)x^{\gamma_j}$ where $a_j+ib_j=0$ if $\gamma_j \notin A_0$.  In order to construct the formal solution we need to solve for the coefficients $a_j(t)+ib_j(t)$ of $x^{\gamma_j}$.  If $\sum_{j=0}^{\infty} \left(a_j\left(t\right)+ib_j\left(t\right)\right) x^{\gamma_j}$ is the positive $x$ formal solution to (\ref{NLS}) then,
\begin{align*}
 \sum_{j=0}^{\infty}& \left[i \dot{a}_j(t)-\dot{b}_j(t)\right] x^{\gamma_j} \\
&=  - \mu \left( \sum_{j=0}^{\infty} \left( \sum_{ {l,m,n}\atop{\gamma_l+\gamma_m+\gamma_n = \gamma_j}}\left[a_l(t)+i b_l(t)\right]\cdot \left[ a_m(t)+i b_m(t) \right] \cdot \left[ a_n(t)-i b_n(t) \right] \right) x^{\gamma_j}\right) \\
&  - \Big(\sum_{j=0}^{\infty}\gamma_j \cdot \left(\gamma_j -1 \right)\cdot \left[a_j(t)+i b_j(t)\right] \cdot x^{\gamma_j-2} \Big)
\end{align*}
from which we deduce that the coefficients $a_j(t), b_j(t)$ satisfy the equations, 
\begin{equation}
\dot{a}_j =  - \mu \sum_{{l,m,n}\atop {\gamma_l+\gamma_m+\gamma_n=\gamma_j}} - a_l\cdot a_m \cdot b_n + a_l\cdot a_n \cdot b_m +a_n\cdot a_m \cdot b_l+b_l\cdot b_m \cdot b_n - \sum_{{p}\atop{\gamma_p-2=\gamma_j}} b_p \cdot \gamma_p \cdot (\gamma_p - 1) \label{a.1}
\end{equation}
\begin{equation}
\dot{b}_j =   \mu \sum_{{l,m,n}\atop {\gamma_l+\gamma_m+\gamma_n=\gamma_j}} a_l\cdot a_m \cdot a_n + a_l\cdot b_n \cdot b_m +b_n\cdot a_m \cdot b_l-b_l\cdot b_m \cdot a_n  + \sum_{{p}\atop{\gamma_p-2=\gamma_j}} a_p \cdot \gamma_p \cdot (\gamma_p - 1) \label{a.2}
\end{equation} 
\indent First we will consider the case when $\gamma_0 < 0$.  Notice first that for $j=0$ the second sums of (\ref{a.1}) and (\ref{a.2}) are nonexistent since $\gamma_0 \geq \gamma_p$ for all $p\geq 0$ and hence there is no $p\geq 0 $ such that $\gamma_p-2 = \gamma_0$.  Also for $j=0$ the first sums are both nonexistent as well because $\gamma_l+\gamma_m+\gamma_n \leq 3\gamma_0 $ and if $\gamma_0 = \gamma_l+\gamma_m+\gamma_n$ then $\gamma_0 \leq 3\gamma_0$ and hence $\gamma_0 \geq 0$ which is a contradiction to our assumption that $\gamma_0 < 0$.  Thus we have $\dot{a}_0=0$ and $\dot{b}_0 = 0$ and hence $a_0(t) = a_0$ is constant and $b_0(t) = b_0$ is constant.  Moreover, for $j\geq 0$ we can see that both sums on the right sides of (\ref{a.1}) and (\ref{a.2}) only contain indices less than $j$.  To see this let us first consider the second sums.  If $\gamma_p-2 = \gamma_j$ then $\gamma_p = \gamma_j+3>\gamma_j$ and hence $p<j$.  For the first sums, if $\gamma_j = \gamma_l+\gamma_m+\gamma_n$ and $l \geq j$ then $0\leq \gamma_j-\gamma_l = \gamma_m+\gamma_n $ so that $\gamma_m+\gamma_n\geq 0$, but $\gamma_m,\gamma_n < 0$, so this is a contradiction, thus $l < j$.  The same argument shows that $m< j$ and $n< j$.  Therefore we may solve for $a_j$ and $b_j$ recursively by integrating the right sides of equations (\ref{a.1}) and (\ref{a.2}) to obtain a polynomials in $t$.  By construction the polynomials will be identically zero for the first few indices until we reach $a_j+ib_j \ne 0$, then it will be a constant $a_j(t)+ib_j(t) = a_j+ib_j$, and for all larger indices $a_j(t)+ib_j(t)$ is polynomial and hence each $a_j(t)+ib_j(t)$ is defined for all $t\in \mathbb{R}$.  \\ 
\indent Now let us assume that $\gamma_0 = 0$.  When $j=0$ the second sums of (\ref{a.1}) and (\ref{a.2}) are again nonexistent for the same reason given above however the first sums are nonzero.  If $\gamma_l+\gamma_m+\gamma_n =\gamma_0 =0$ then $\gamma_l = \gamma_m = \gamma_n = 0$.  Therefore $a_0(t)$ and $b_0(t)$ satisfy the equations 
\begin{eqnarray*}
\dot{a}_0 &=& -\mu \left(a_0^2 b_0+b_0^3\right) \\
\dot{b}_0 &=& \mu\left(a_0^3+a_0 b_0^2\right)
\end{eqnarray*} 
which can easily be solved to yield the solutions 
\begin{eqnarray*}
a_0(t) & = & a_0(0)\cos(c\mu t) + b_0(0) \sin(c\mu t) \\
b_0(t) & = & -a_0(0)\sin(c\mu t) + b_0(0) \cos(c\mu t)
\end{eqnarray*}
where $c=a_0(0)^2+b_0(0)^2$ and hence $a_0(t)$ and $b_0(t)$ are continuous and nonvanishing for $t\in \mathbb{R}$.  \\
\indent Now we'll consider the construction of functions $a_j(t),b_j(t)$ for $j \geq 1$.  When $j \geq 1$ the second sums of (\ref{a.1}) and (\ref{a.2}) always consist of indices less than $j$ however for the first sums all terms involve indices which are at most $j$.  To see this, suppose that $\gamma_j = \gamma_l+\gamma_m+\gamma_n $ and suppose that $l>j$.  Then $0 < \gamma_j-\gamma_l = \gamma_m+\gamma_n$ and hence $\gamma_m+\gamma_n > 0$ which is a contradiction since $\gamma_m,\gamma_n \leq 0$.  Therefore we conclude that if $\gamma_j = \gamma_l+\gamma_m+\gamma_n$  and $j\geq 1$ then $l,m,n \leq j$.  \\
\indent Furthermore, if $l=j$ then $\gamma_m+\gamma_n = 0$ which implies that $m=n=0$.  Similar statements hold when $m=j$ and when $n=j$.  Therefore when $j \geq 1$ the coefficients $a_j(t),b_j(t)$ satisfies the equations 
\begin{eqnarray*}
\dot{a}_j & = &  -\mu \left( 2a_0 b_0 a_j + a_0^2 b_j + 3 b_0^2 b_j\right)+ P_j(a_0,b_0\ldots,a_{j-1},b_{j-1}) \\
\dot{b}_j & = & \mu \left( 3a_j a_0^2 +a_j b_0^2 + 2a_0 b_0 b_j \right) + Q_j(a_0,b_0\ldots,a_{j-1},b_{j-1})
\end{eqnarray*}
where $P_j$ and $Q_j$ are polynomials and hence the functions $a_j(t)$ and $b_j(t)$ exist and are continuous for $t\in \mathbb{R}$.     
This concludes the construction of the formal solution $\sum_{j=0} ^{\infty} \left[a_j ^{\pm}(t)+ib_j^{\pm}(t)\right] x^{\gamma_j}$.  \\
\indent  Now let $f(x,t)$ denote any smooth function which is asymptotic to the formal solution \\ $\sum_{j=0} ^{\infty} \left[a_j ^{\pm}(t)+ib_j^{\pm}(t)\right] x^{\gamma_j}$ (by proposition 3.5 in \cite{Shub} there exists such a function).  By plugging in $f(x,t)$ to (\ref{NLS}) we will now show that one obtains a function $g(x,t)\in S^{-\infty}(\mathbb{R}\times \mathbb{R} \to \mathbb{C})$.  Suppose $\left|x\right|\geq 1$ and $J\subset \mathbb{R}$ is a compact subset and  $l,m,N\geq 0$ are integers.  Let $S_N(f) = \sum_{j=0} ^{N} \left[a_j ^{\pm}(t)+ib_j^{\pm}(t)\right] x^{\gamma_j}$ (i.e. if $x>1$ then $S_N(f) = \sum_{j=0} ^{N} \left[a_j ^{+}(t)+ib_j^{+}(t)\right] (+x)^{\gamma_j}$ and when $x<-1$ we have $S_N(f) = \sum_{j=0} ^{N} \left[a_j ^{-}(t)+ib_j^{-}(t)\right] (-x)^{\gamma_j}$ ).  Then we have

\begin{eqnarray*}
\partial^l_t\partial^m_x \left( i f_t+ f_{xx}+\mu \overline{f} f^2  \right) & = & i \partial^l_t\partial^m_x \left[S_N(f) +\left(f-S_N(f)\right) \right]_t + \partial^l_t\partial^m_x\left[S_N(f) +\left(f-S_N(f)\right) \right]_{xx} \\
                     &    &  + \mu \partial^l_t\partial^m_x \left[\overline{\left(S_N(f) +\left(f-S_N(f)\right) \right)} \left( S_N(f) +\left(f-S_N(f)\right) \right)^2\right]
\end{eqnarray*}

\indent After expanding the right side we will obtain the expression \\ $\partial^l_t\partial^m_x \left[ i S_N(f)_t+ S_N(f)_{xx}+\mu \overline{S_N(f)} S_N(f)^2 \right]$ and terms which are products constants times \\ $\partial^p_t \partial^q_x \left(f-S_N(f) \right)^r$ and $ \partial^s_t \partial^k_x S_N(f)^n$ for some $p,q,r,s,n \in \mathbb{N}$.  Since the coefficients $a_j(t), b_j(t)$ satisfy the equations (\ref{a.1}) and (\ref{a.2}) and since $f$ is asymptotic to $\sum_{j=0} ^{\infty} \left[a_j ^{\pm}(t)+ib_j^{\pm}(t)\right] x^{\gamma_j}$ we have for any $M\in \mathbb{N}$ there exists a sufficiently large $N$ such that these terms are bounded by $C_{M,J,l,m}\left|x \right|^{-M}$.\hspace{1cm}$\square$ \\





\end{document}